\documentclass[a4]{article}

\usepackage{amsmath}
\usepackage{amssymb}
\usepackage{amsthm}
\usepackage{amsxtra}
\usepackage{pstcol}
\usepackage{graphicx}
\usepackage{ifthen}

\def\de{Definition}
\def\pp{Proposition}
\def\th{Theorem}

\def\lm{Lemma}
\def\co{Corollary}
\def\ex{Example}
\def\exs{Examples}
\def\re{Remark}

\newcommand{\comment}[1]{}

\newcommand{\is}{\ensuremath{\cap}}

\renewcommand{\vert}{\ensuremath{\,|\,}}

\renewcommand{\deg}[2][{}]{\ifthenelse{\equal{#1}{}}
				{\normalfont\ensuremath{\textrm{deg}(#2)}}
				{\normalfont\ensuremath{\textrm{deg}_{\mathcal #1}(#2)}}}

\renewcommand{\det}[1]{\normalfont\ensuremath{\textrm{det}(#1)}}
\newcommand{\ordP}[1]{\normalfont\ensuremath{\textrm{ord$_P$}(#1)}}
\newcommand{\Spec}[1]{\normalfont\ensuremath{\textrm{Spec}(#1)}}
\newcommand{\Max}[1]{\normalfont\ensuremath{\textrm{Max}(#1)}}
\newcommand{\Hess}[1]{\normalfont\ensuremath{\textrm{Hess}_#1}}

\newcommand{\Sing}[1]{\normalfont\ensuremath{\textrm{Sing}(#1)}}

\newcommand{\nfrac}[2]
{\raisebox{0.5ex}{\footnotesize \ensuremath{#1}}\hspace{-0.1ex}\raisebox{0.1ex}{/}\hspace{-0.2ex}\raisebox{-0.4ex}{\footnotesize \ensuremath{#2}}}

\newcommand{\demph}[1]{\emph{#1}\Index{#1}}

\newcommand{\Nat}{\ensuremath{\mathbb N}}

\newcommand{\ifdef}[2]{\ifthenelse{\boolean{#1}}{#2}{}}

\def\hA{\section}

\def\hAA{\section*}

\def\Index#1{}

\newcommand{\headA}[1]{\hA{#1}\Index{#1}}

\newcommand{\headAA}[1]{\hAA{#1}\Index{#1}\addcontentsline{toc}{section}{#1}
}

\theoremstyle{definition}
\newtheorem{Definition}{\de}[section]
\newtheorem{Example}[Definition]{\ex}
\newtheorem{Examples}[Definition]{\exs}
\newtheorem{Remark}[Definition]{\re}

\theoremstyle{plain}
\newtheorem{Proposition}[Definition]{\pp}
\newtheorem{Lemma}[Definition]{\lm}
\newtheorem{Corollary}[Definition]{\co}
\newtheorem{Theorem}[Definition]{\th}

\newcommand{\defn}[2][{}]{\ifthenelse{\equal{#1}{}}
				{\begin{Definition}#2\end{Definition}}
				{\begin{Definition}[#1]#2\end{Definition}}}
\newcommand{\prop}[2][{}]{\ifthenelse{\equal{#1}{}}
				{\begin{Proposition}#2\end{Proposition}}
				{\begin{Proposition}[#1]#2\end{Proposition}}}
\newcommand{\lma}[2][{}]{\ifthenelse{\equal{#1}{}}
				{\begin{Lemma}#2\end{Lemma}}
				{\begin{Lemma}[#1]#2\end{Lemma}}}
\newcommand{\cor}[2][{}]{\ifthenelse{\equal{#1}{}}
				{\begin{Corollary}#2\end{Corollary}}
				{\begin{Corollary}[#1]#2\end{Corollary}}}
\newcommand{\thm}[2][{}]{\ifthenelse{\equal{#1}{}}
				{\begin{Theorem}#2\end{Theorem}}
				{\begin{Theorem}[#1]#2\end{Theorem}}}
\newcommand{\rem}[2][{}]{\ifthenelse{\equal{#1}{}}
				{\begin{Remark}#2\end{Remark}}
				{\begin{Remark}[#1]#2\end{Remark}}}
\newcommand{\expl}[2][{}]{\ifthenelse{\equal{#1}{}}
				{\begin{Example}#2\end{Example}}
				{\begin{Example}[#1]#2\end{Example}}}
\newcommand{\expls}[2][{}]{\ifthenelse{\equal{#1}{}}
				{\begin{Examples}#2\end{Examples}}
				{\begin{Examples}[#1]#2\end{Examples}}}
\newcommand{\prf}[2][{}]{\ifthenelse{\equal{#1}{}}
				{\begin{proof}#2\end{proof}}
				{\begin{proof}[#1]#2\end{proof}}}

\begin{document}

\begin{center}
{\huge\bf Weierstra\ss\ semigroups and nodal curves of type $p,q$}
\end{center}
\begin{center}
	{H.\ Knebl}
	\footnote{Fakult\"at Informatik, Georg-Simon-Ohm-Hochschule N\"urnberg, Kesslerplatz 12, 90489 N\"urnberg, Germany},
	{E.\ Kunz $^2$ and R.\ Waldi}
	\footnote{Fakult\"at f\"ur Mathematik, Universit\"at Regensburg, Universit\"atsstrasse 31, 93053 Regensburg, Germany}
\end{center}
\vspace{0.5cm}
\begin{abstract}
We study plane curves of type $p,q$ having only nodes as singularities. 
Every Weierstra\ss\ semigroup is the Weierstra\ss\ semigroup of such a curve at its place at infinity for properly chosen $p,q$. 
We construct plane curves of type $p,q$ with explicitly given nodes and determine their Weierstra\ss\ semigroups. 
Many such semigroups are found.
\end{abstract}

\headAA{Introduction}\label{Introduction}

Let $\mathcal R$ be a smooth projective algebraic curve of genus $g$ defined over an algebraically closed
field $K$ of characteristic 0. For a closed point $P\in \mathcal R$ consider the meromorphic functions on
$\mathcal R$ which have a pole at most at $P$, hence are regular everywhere else. By the
Weierstra\ss\ Gap Theorem the set $H(P)$ of pole orders of such functions is a numerical semigroup of genus $g$,
i.e.\ it is closed under addition and the {\it set of gaps} $\mathbb N\setminus H(P)$ consists of exactly $g$
integers. Moreover for all but finitely many $P\in \mathcal R$ the gaps are $1,\dots,g$.
The semigroup $H(P)$ is called the {\it Weierstra\ss\ semigroup} of $\mathcal R$ at $P$.

Hurwitz [Hu] asked in 1893: Which numerical semigroups do occur as Weierstra\ss\ semigroups?

Since then many classes of Weierstra\ss\ semigroups have been specified. For example all numerical semigroups of
genus $\le 8$ are Weierstra\ss\ semigroups ([Ko1], [Ko-O]), as are complete intersection semigroups ([Pi]).
The first correct example of a numerical semigroup which is not a Weierstra\ss\ semigroup was found by
Buchweitz [B] in 1980. Such semigroups are now called {\it Buchweitz semigroups}.
K.O.\ St\"ohr and F. Torres ([T], Scholium 3.5) have shown that for any $g\ge 100$ {\it symmetric}
Buchweitz semigroups of genus $g$ exist, and in [Ko2] many Buchweitz semigroups are constructed.

In this paper we investigate Weierstra\ss\ semigroups using plane curves of type $p,q$.
For relatively prime integers $p,q$ with $1<p<q$ these are the curves $C$ in $\mathbb A^2(K)$
defined by a {\it Weierstra\ss\ equation of type $p,q$}
$$
	C: Y^p + aX^q + \sum_{\nu p + \mu q < pq} a_{\nu\mu}X^{\nu}Y^{\mu} = 0 \qquad (a_{\nu\mu}\in K, a \in K\setminus\{0\}).
$$
Such curves are irreducible  and have only one place $P$ at infinity which is a point on the
normalization $\mathcal R$ of the projective closure of $C$. We call $H(P)$ also the
{\it Weierstra\ss\ semigroup of $C$}. It contains the semigroup $H_{pq}:=\left\langle p,q \right\rangle$ generated by $p$ 
and $q$ as a subsemigroup. 
The method to study compact Riemann surfaces by investigating their plane models of type $p,q$
can be attributed to Weierstra\ss\ (see the remark on page 544 of [HL]).

The gaps of $H(P)$ can be described by the holomorphic differentials $\omega$ on $\mathcal R$:
$$
	\mathbb N\setminus H(P)=\{\ordP {\omega} + 1 \vert \omega\ \text{holomorphic on } \mathcal R\}
$$
where $\textrm{ord}_P$ is the normed discrete valuation at the point $P$. The holomorphic differentials on $\mathcal R$ 
are closely connected to the adjoint curves of $C$. If $C$ has only nodes as singularities, i.e. singularities of multiplicity 2 
with two different tangents, these are the curves passing through all of the nodes. If the nodes are explicitly known the 
holomorphic differentials on $\mathcal R$ can be computed (Proposition \ref{PropLinEq}), and so can be the Weierstra\ss\ semigroups $H(P)$. 
For this reason we are mainly interested in curves of type $p,q$ having only nodes as singularities (at finite distance). 
We call them \demph{nodal curves of type} $p,q$. 
It can be shown that any Weierstra\ss\ semigroup is the Weierstra\ss\ semigroup of a nodal curve of type $p,q$ 
for properly chosen integers $p$ and $q$ (Theorem \ref{NodalModels}).

Given a nodal curve $C$ of type $p,q$ with explicitly known nodes there is a construction to eliminate some of its nodes. 
More precisely this means that for certain subsets $\{P_1,\dots,P_l\}$ of the set of nodes of $C$ another curve $C'$ of type $p,q$ 
can be found having $\{P_1,\dots,P_l\}$ as its set of nodes (Corollary \ref{CorNodeElim}). By elimination of nodes a great number of 
Weierstra\ss\ semigroups and explicitly given nodal curves of type $p,q$ having these Weierstra\ss\ semigroups can be 
specified, see Section \ref{ConstructionOfWeierstrassSemigroups}.

In our considerations the numerical semigroups containing $p$ and $q$ play a crucial role. We begin with a geometric 
illustration of such semigroups.

\headA{An illustration of numerical semigroups}\label{AnIllustration}

Let $p,q\in \mathbb N$ be relatively prime integers with $1<p<q$, and let $H_{pq}=$ $\left\langle p,q \right\rangle$ be the numerical
semigroup generated by $p$ and $q$. It is a symmetric semigroup with the conductor $c:=(p-1)(q-1)$ and $d:= \nfrac{c}{2}$ gaps 
$\gamma_1<\dots <\gamma_d$ where $\gamma_i=c-1-(a_ip+b_iq)$ with
uniquely determined $a_i,b_i\in \mathbb N$. Therefore the gaps are in one-to-one correspondence with the $d$
lattice points $(a,b)\in \mathbb N^2$ below the line
$$g_0: pX+qY=c-1.$$
Let $\Delta _0$ be the triangle which is determined by $g_0$ and the coordinate axes. On each parallel
$$g_i: pX+qY=c-1-i\qquad (i=1,\dots,c-1)$$
there is at most one lattice point, and it corresponds to the gap $i$ of $H_{pq}$. Let $\Delta_i$ be the triangle
determined by $g_i$ and the coordinate axes.

\def\PicScale{0.35}
\begin{center}
	\includegraphics[clip,angle=-90,scale=\PicScale]{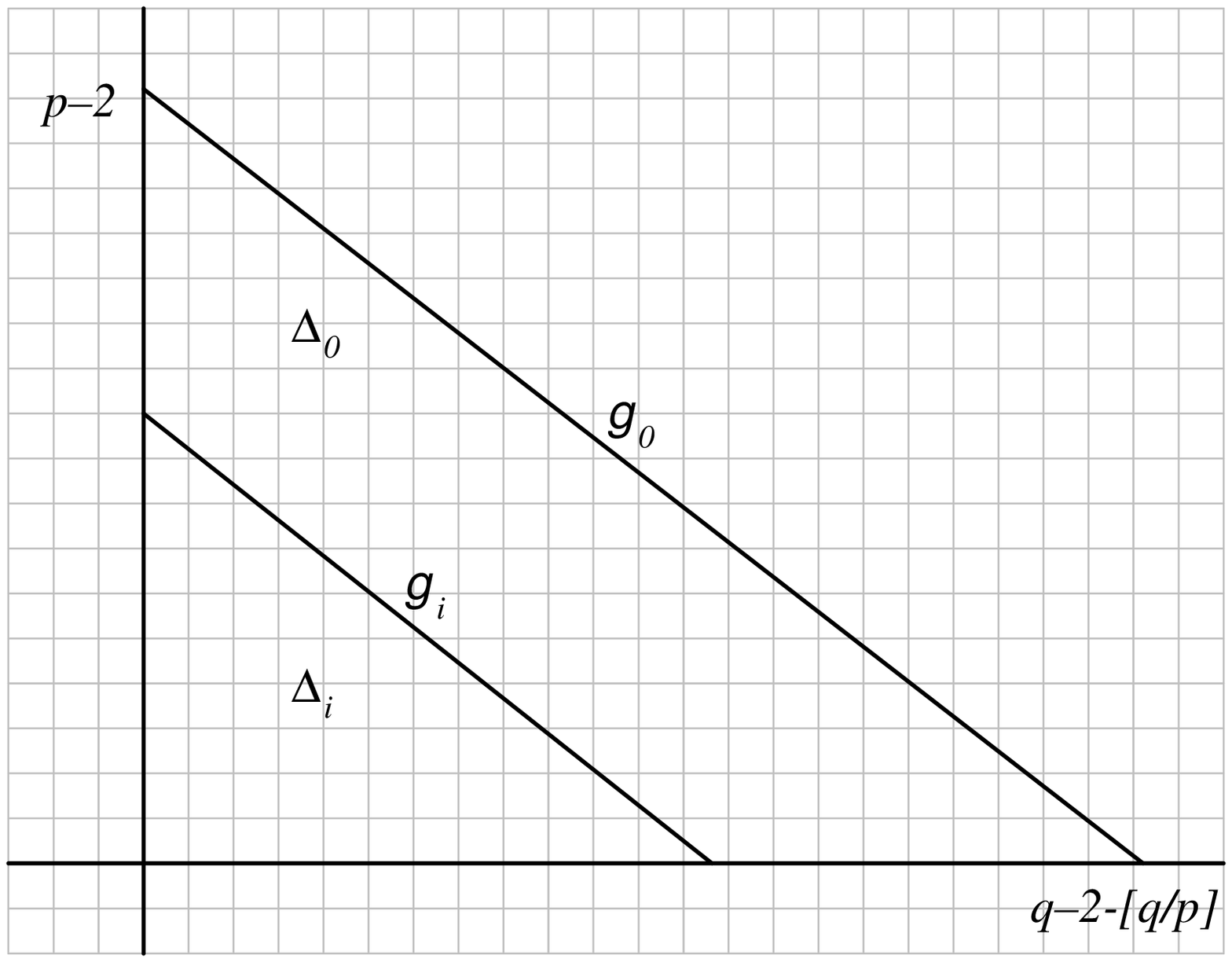}\hbox{}
\end{center}

A numerical semigroup $H$ with $p,q\in H$ arises from $H_{pq}$ by closing some of its gaps. If
$(a,b)\in \Delta_1$ corresponds to such a gap, then also the lattice points of the rectangle $R_{a,b}$ with the
corners $(0,0),(a,0),(a,b),(0,b)$ belong to gaps which are closed in $H$. Therefore these gaps correspond
to a set $L_H\subset \Delta_1$ of lattice points determined by a lattice path starting on the $Y$-axis
and ending on the $X$-axis with downward and right steps, see the following figure, where $P_0$ and $P_1$ 
have the same $Y$-coordinates if the first step in the lattice path is a right step, and $P_{m - 1}$ and $P_m$ may have
the same $X$-coordinates.

\begin{center}
	\includegraphics[clip,angle=-90,scale=\PicScale]{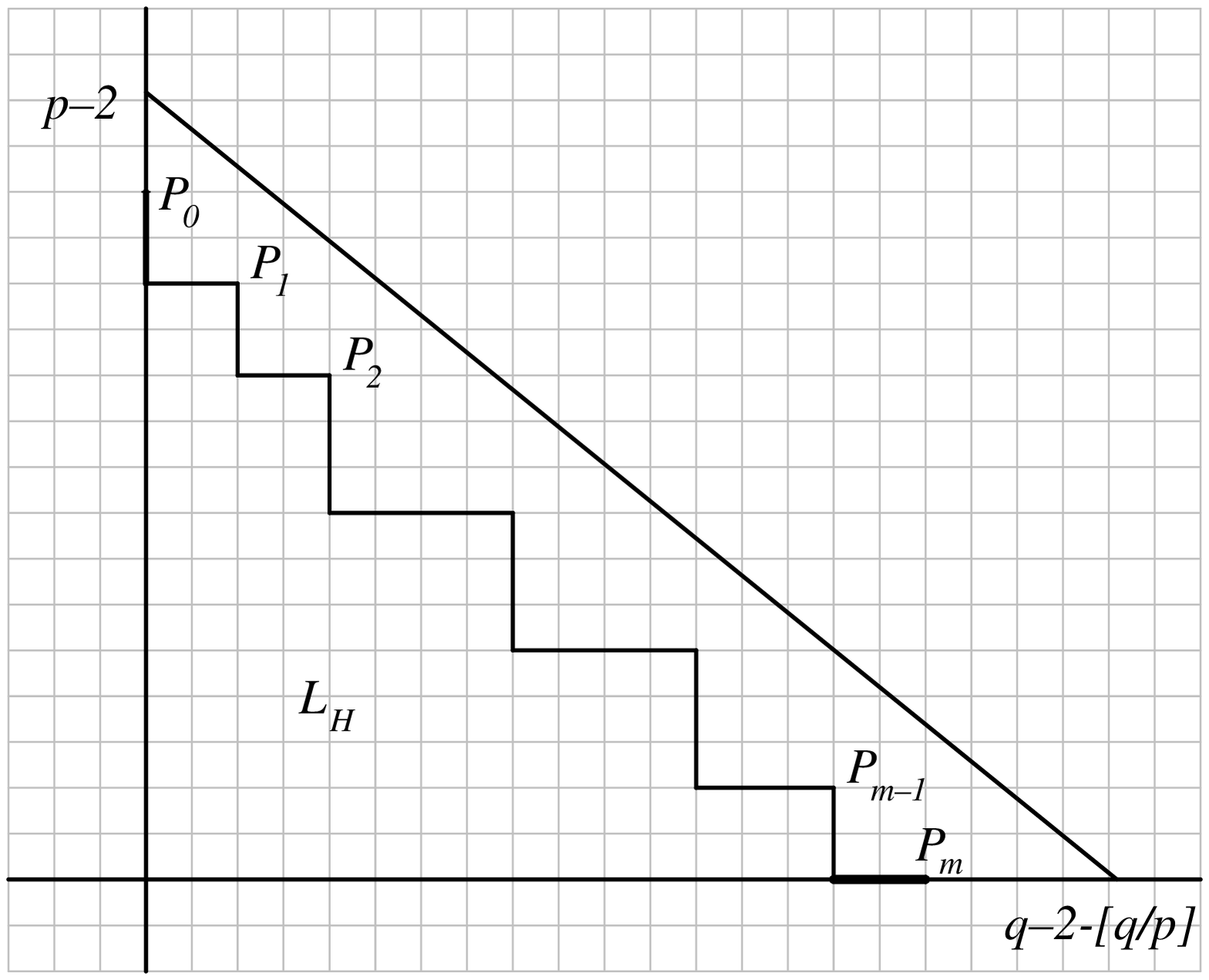}\hbox{}
\end{center}

If $P_i=(\alpha_i,\beta_i)\ (i=0,\dots,m)$, then
$$
	\alpha_0=0 < \alpha_1 < \dots  < \alpha_{m-1}\le \alpha_m \le q-2-\left\lfloor\frac{q}{p}\right\rfloor
$$
$$
	p-2 \ge \beta_0 \ge \beta_1 > \dots > \beta_{m-1} > \beta_m=0.
$$
The numbers $c-1-(\alpha_ip+\beta_iq)\ (i=0,\dots,m)$ together with $p$ and $q$ form a system of generators of $H$.

We write $L=(P_0,\dots,P_m)$ for such subsets of lattice points in $\Delta_1$. Not every $L=(P_0,\dots,P_m)$ belongs to a
numerical semigroup. For two gaps $\gamma,\gamma'$ of $H_{pq}$ which are to be closed in
$H$ also $\gamma+ \gamma'$ must be such a gap or an element of $H_{pq}$.

\expls{\label{ExIll_1}

\noindent
a) Let $P=(a,b)\in g_1,\ a\ne 0,b\ne 0$. Then the set $L=((0,b),P,(a,0))$ does not correspond to a numerical semigroup, since with 1 all 
gaps of $H_{pq}$ are closed in $H$.

\noindent
b) If $L$ consists of all lattice points in $\Delta_i$, then $L=L_H$ where $H$ is the semigroup obtained from $H_{pq}$ by closing all its 
gaps which are $\ge i$.
}

\lma{
Let $\gamma$ and $\gamma'$ be gaps of $H_{pq}$ belonging to $(a,b)\in \Delta_1$ and $(a',b')\in \Delta_1$ respectively. 
Then $\gamma+\gamma'$ is a gap of $H_{pq}$ if and only if $a+a'\ge q-1$ or $b+b'\ge p-1$. The gap $\gamma+\gamma'$ corresponds to 
$Q:=(a,b)+(a',b')+(1-q,1)$ if $a+a'\ge q-1$, and to $P:=(a,b)+(a',b')+(1,1-p)$ if $b+b'\ge p-1$.
}

\prf{
We have
$$\gamma+\gamma'=c-1+pq-p-q-[(a+a')p+(b+b')q]\ge 0$$
Write

$$\gamma+\gamma'=c-1-[(a+a'+1)p+(b+b'-p+1)q]\leqno(1)$$
$$\gamma+\gamma'=c-1-[(a+a'-q+1)p+(b+b'+1)q].\leqno(2)$$

If $a+a'\ge q-1$ or $b+b'\ge p-1$, then $\gamma+\gamma'$ is a gap of $H_{pq}$, and it belongs to $Q\in \Delta_1$ or $P\in \Delta_1$ respectively.

Conversely let $\gamma+\gamma'$ be a gap of $H_{pq}$ belonging to $(a'\!',b'\!')\in \Delta_1$. Then

$$(a+a'-q+1)p+(b+b'+1)q=a'\!'p+b'\!'q.\leqno(3)$$

This implies
$$q\vert a'\!'-a-a'-1\ \text{and}\ p\vert b'\!'-b-b'-1.$$
As $a,a',a'\!'<q-1$ and $b,b',b'\!'<p-1$ we obtain

$$a'\!'-a-a'-1\in \{0,-q\}\ \text{and}\ b'\!'-b-b'-1\in \{0,-p\}.$$
But $(a'\!',b'\!')\ne (a+a'+1,b+b'+1)$ by (3). Hence
$$a+a'+1-q=a'\!'\ge 0\ \text{or}\ b+b'+1-p=b'\!'\ge 0$$
and $a+a'\ge q-1$ or $b+b'\ge p-1$.
}

We obtain

\prop{\label{ClosedCond}
$L=(P_0,\dots,P_m)$ belongs to a numerical semigroup if and only if $L$ is closed
under the operations
$$(a,b),(a',b')\mapsto
\begin{cases}
    (a,b)+(a',b')+(1-q,1), \text{if}\ a+a'\ge q-1,\\
    (a,b)+(a',b')+(1,1-p), \text{if}\ b+b'\ge p-1.
\end{cases}$$
If $P_i=(\alpha_i,\beta_i)\ (i=0,\dots,m)$, then it suffices to check the conditions for the $(\alpha_i,\beta_i)$ and
$(\alpha_j,\beta_j)$ with $\alpha_i+\alpha_j\ge q-1$ or $\beta_i+\beta_j\ge p-1$.
}

\hspace{-0.5cm}
\begin{minipage}[t]{12 cm}
	\begin{minipage}[c]{10em}
	\includegraphics[clip,angle=-90,scale=\PicScale]{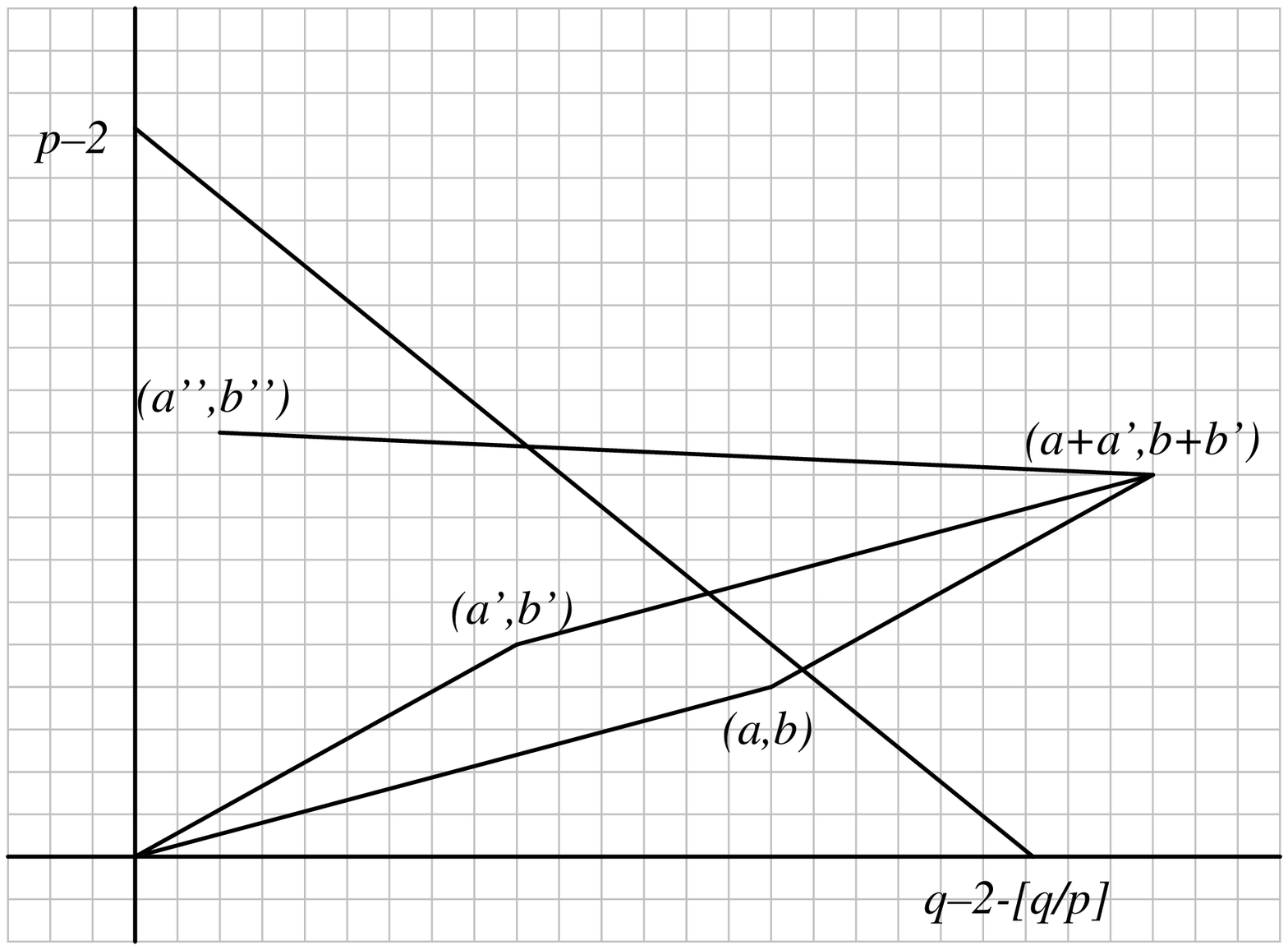}
	\end{minipage} \hfill
	\begin{minipage}[c]{10em}
	\includegraphics[clip,angle=-90,scale=\PicScale]{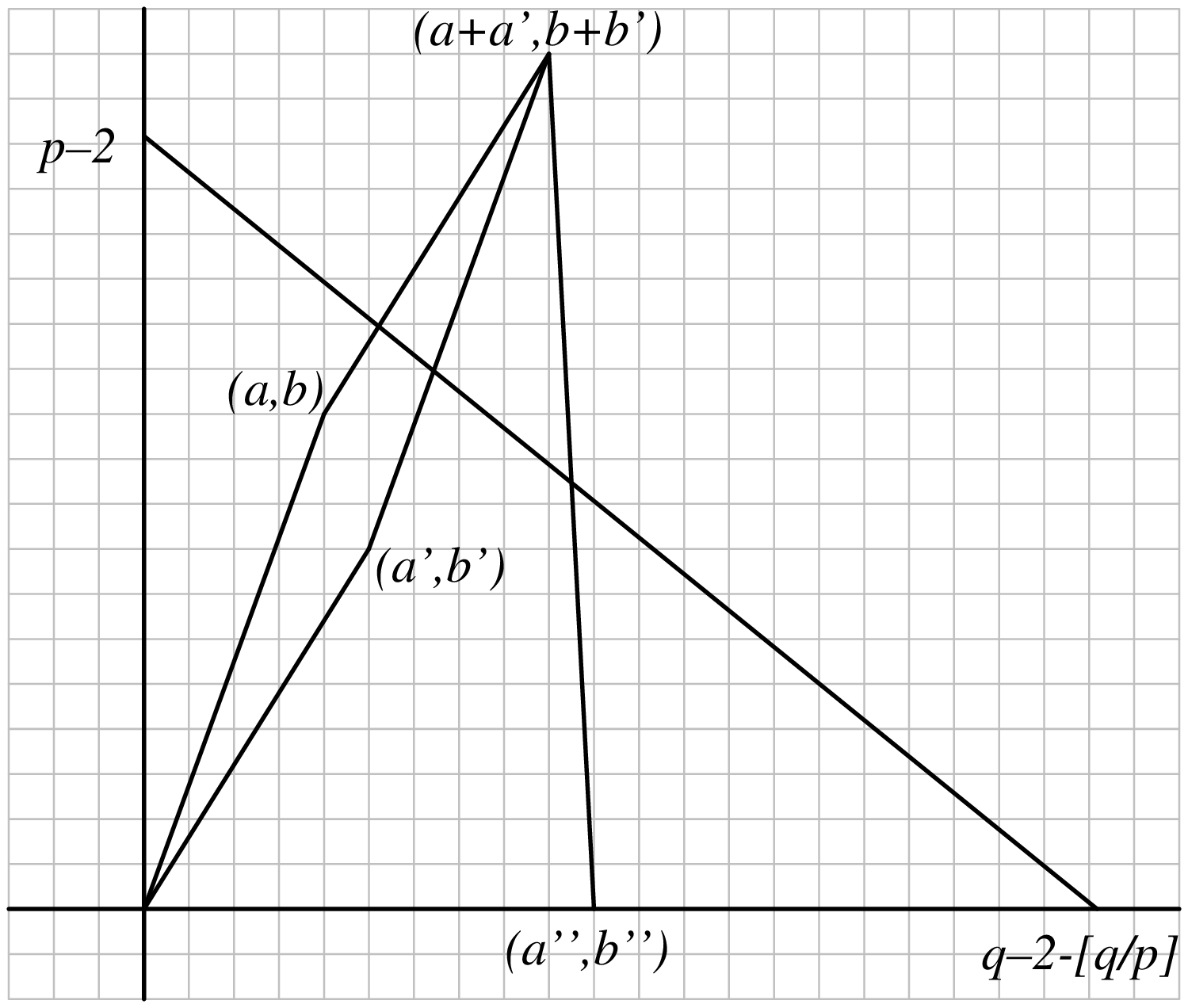}
	\end{minipage}
\end{minipage}

\expls{\label{ExSemigroups}
\noindent
a) A rectangle $R_{a,b}$ as above defines a numerical semigroup $H$ if and only if $a<\nfrac{1}{2}(q-1), b<\nfrac{1}{2}(p-1)$. 
Then $H=  \left\langle p,q,c-1-(ap+bq) \right\rangle =: H_{(a,b)}$. If the conditions are satisfied and $L$ is defined by a lattice path inside the 
rectangle, then $L=L_H$ with a numerical semigroup $H\subset H_{(a,b)}$, since the closedness condition of the proposition 
is clearly fulfilled.

\begin{center}
	\includegraphics[clip,angle=-90,scale=\PicScale]{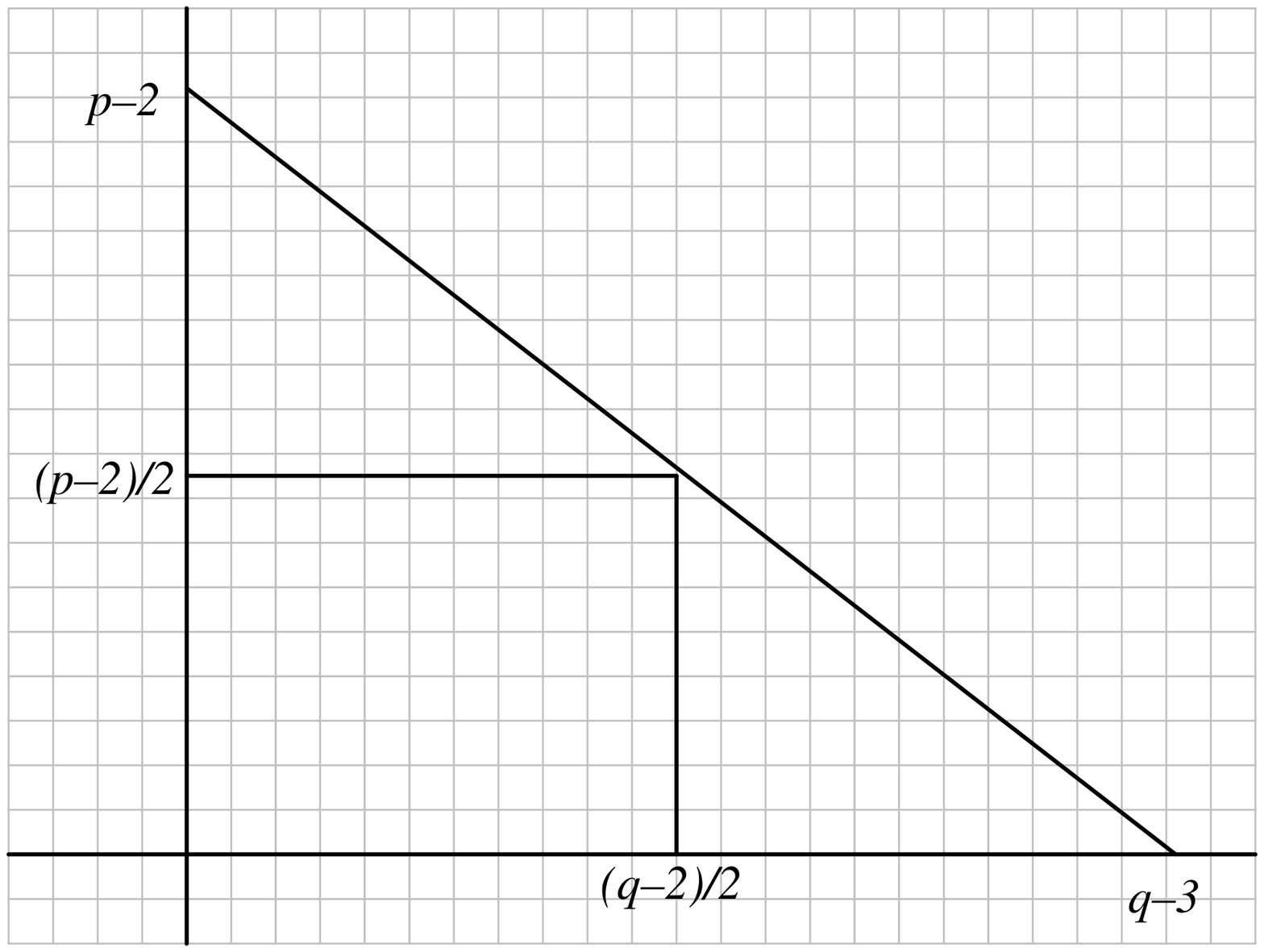}\hbox{}
\end{center}

\noindent
b) Let $\Delta_k\subset L\subset \Delta_i$ with $2i\ge k$ and $i\ge 1$. Then $L=L_H$ with a
        numerical semigroup $H$.
In fact, let $(a,b),(a',b')\in L$ and $a+a'\ge q-1$. Then from $ap+bq\le c-1-i$ and
        $a'p+b'q\le c-1-i$ we obtain
        $$(a+a'+1-q)p+(b+b'+1)q\le c-1-2i\le c-1-k$$
        hence $(a,b)+(a',b')+(1-q,1)\in \Delta_k\subset L$. Similarly if $b+b'\ge p-1$.

\noindent
c) For $r\in \mathbb N$ with $0\le r \le p-2$ let $L$ be the triangle with the corners $(0,0),(r,0),(0,r)$. 
Then $L=L_H$ with a numerical semigroup $H$. In fact, if $(a,b), (a',b')\in L$, then $a+b\le r,a'+b'\le r$. 
If $a+a'\ge q-1$, then we have for $(a'\!',b'\!'):=(a,b)+(a',b')+(1-q,1)$ that
$$
	a'\!'+b'\!'=a+a'-q+b+b'+2\le2r-(q-2)\le r.
$$ 
Similarly, if $b+b'\ge p-1$.

\noindent
d) Let $q=p+1$ and $H_p:=H_{pq}$. In this case the gaps of $H_p$ are in one-to-one correspondence with the
$\binom{p}{2}=\nfrac{c}{2}$ lattice points of the triangle $\Delta$ with the corners $(0,0),(p-2,0),(0,p-2)$,
that is with the $(a,b)\in \mathbb N^2$ with $a+b < p-1$. Let $s_i$ be the section from $(0,p-2-i)$ to
$(p-2-i,0)\ (i=0,\dots,p-2)$. The lattice points on this section correspond to the gaps $ip+i+1,\dots,ip+p-1$
of $H_p$. Let $\Delta_i'$ be the triangle with the corners $(0,0),(p-2-i,0),(0,p-2-i)$ and let
$\Delta_k'\subset L\subset \Delta_i'$ with $2i\ge k$, then $L$ defines a numerical semigroup because the
conditions of Proposition \ref{ClosedCond} are fulfilled. Especially this is true if we add to the lattice points of
$\Delta_k'\ (k\ge 2)$ some lattice points of $s_{k-1}$.
}

It is a much more difficult problem to decide which $L$ define Weierstra\ss\ semigroups. Later in the paper we will describe 
classes of such $L$ for which this is the case.

On the other hand, take for instance $p\ge 9,q=p+1$. By \ref{ExSemigroups}d) the numbers $1,\dots,p-1,2p-7,2p-5,2p-2,2p-1$ are the gaps of a 
numerical semigroup $H$ of genus $g=p+3$ with $p,p+1\in H$. Similarly for the numbers $1,\dots,p-1,2p-7,2p-6,2p-3,2p-1$.
Let $l_2(H)$ denote the number of sums of two gaps of $H$. If $p\ge 13$, then in both cases $l_2(H)=3p+7>3g-3$,
hence the necessary condition for Weierstra\ss\ semigroups given by Buchweitz
$$
	l_2(H)\le 3g-3\leqno(4)
$$
is hurt. For larger $p$ there is an increasing number of possibilities to construct Buchweitz semigroups.
The work by Komeda [Ko1] gives all semigroups of genus 16 and 17 for which (4) is hurt and contains a list
with the numbers of all semigroups of genus $\le 37$ and with the numbers of those for which (4) is hurt.
For $p=13$ the two examples above are the semigroups of genus 16 for which (4) is hurt, the first one is
the example found by Buchweitz [B] of a numerical semigroup which is not a Weierstra\ss\ semigroup.

\headA{Plane algebraic curves of type $p,q$}\label{PlaneAlgebraic}

For relatively prime numbers $p,q\in \mathbb N$ with $1<p<q$
let $C: F=0$ be a curve of type $p,q$ in $\mathbb A^2(K)$ with 
defining Weierstra\ss\ polynomial 
$$
	F := Y^p+aX^q+\sum_{\nu p+\mu q<pq}a_{\nu\mu}X^{\nu}Y^{\mu}.
$$

\prop{\label{OnePlaceInf}
$C$ is irreducible and has only one place $P$ at infinity, hence also only one 
point at infinity. If $x,y\in K[C]=K[X,Y]/(F)$ are the residue classes of $X,Y$, then
$$
	\ordP x = -p,\ \ordP y = -q.
$$
}

\prf{
([Kn], Satz 1)
}

\prop{\label{SingDeg}
Let $\bar C$ be the projective closure of $C$ and $Q$ its point at infinity. At $Q$ 
the curve $\bar C$ has multiplicity $q-p$ and singularity degree
$$
	\delta_Q=\frac{1}{2}(q-p-1)(q-1).
$$
}
\prf{
([Kn], Satz 2)
}
\cor{\label{GenusFormula}
If $C$ has $l$ nodes and no other singularities (at finite distance), then the normalization $\mathcal R$ of 
$\bar C$ has genus
$$
	g = \frac{1}{2}(p-1)(q-1)-l.
$$
}

\prf{
This follows from the genus formula
$$
	g=\frac{1}{2}(q-1)(q-2)-\sum_{R\in \tilde C(K)}\delta_R
$$
for plane projective curves of degree $q$ (see e.g.\ [Ku], Theorem 14.7) and the fact that $\delta_R=1$ for a node $R$ and 
$\delta_Q=\nfrac{1}{2}(q-p-1)(q-1)$ for the point $Q$ at infinity.
}

By Proposition \ref{OnePlaceInf} the Weierstra\ss\ semigroup $H(P)$ of a curve $C$ of type $p,q$ contains $p$ and $q$, 
hence we are in the situation of Section \ref{AnIllustration}. If $C$ is a nodal curve, then by Corollary \ref{GenusFormula} 
the number of its nodes is at most $\nfrac{1}{2}(p-1)(q-1)$.

\expls[Lissajous curves]{\label{LissajousCurves}
\noindent
The Chebyshev polynomials $T_n(X)$ of the first kind are defined by
$$
	T_0(X):=1,\quad T_1(X):=X,\quad T_n(X):=2XT_{n-1}(X)-T_{n-2}(X)\ \text{for}\ n\ge 2.
$$
They have degree $n$ and leading coefficient $2^{n-1}$. If $p$ and $q$ are relatively prime integers with $1<p<q$, then 
the {\it Lissajous curve $C_L$ of type $p,q$} is the plane curve defined over $\mathbb C$ by the polynomial 
$$
	L_{p,q}(X,Y):=T_p(Y)-T_q(X).
$$
Up to a factor $2^{p-1}$ it is a Weierstra\ss\ polynomial of type $p,q$.
As $T_n(T_m(X))=T_{nm}(X)$ the curve $C_L$ has the parametric representation $x=T_p(t),\ y=T_q(t)$.
It has the nodes
$$(x_k,y_l):= \left(\cos\frac{k\pi}{q},\cos\frac{l\pi}{p}\right)\ (k=1,\dots,q-1, l=1,\dots,p-1, k\equiv l\ \text{mod}\ 2).$$
Their number is $d:=\nfrac{1}{2}(p-1)(q-1)$.

In fact, we have $T'_p(Y)=pU_{p-1}(Y), T'_q(X)=qU_{q-1}(X)$ with the Chebyshev polynomials $U_{p-1}, U_{q-1}$ of the second kind. 
Here $U_{q-1}$ has simple zeros $x_k\ (k=1,\dots,q-1)$ and $U_{p-1}$ has simple zeros $y_l\ (l=1,\dots, p-1)$ so that the singularities of 
$C_L$ are found among the $(x_k,y_l)\ (k=1,\dots,q-1, l=1,\dots,p-1)$. The Hesse determinant of $L_{p,q}$ does not vanish at the $(x_k,y_l)$ 
since the $x_k,y_l$ are simple zeros of $T'_q$ resp. $T'_p$. Therefore the singularities are nodes. Moreover it is known that
$$T_q\left(\cos\frac{k\pi}{q}\right)=\cos(k\pi)=(-1)^k,\quad T_p\left(\cos\frac{l\pi}{p}\right)=\cos(l\pi)=(-1)^l$$
and it follows that $L_{p,q}(x_k,y_l)=T_p(y_l)-T_q(x_k)=0$ if and only if $k \equiv l \mod 2$.\\

The real points of the Lissajous curves defined by
\begin{eqnarray*}
	L_{4,7} &=& 8Y^4-8Y^2+1-64X^7+112X^5-56X^3+7X, \\
	L_{5,7} &=& 16Y^5-20Y^3+5Y-64X^7+112X^5-56X^3+7X \hbox{ and} \\
	L_{5,8} &=& 16Y^5-20Y^3+5Y-128X^8+256X^6-160X^4+32X^2-1
\end{eqnarray*}

are sketched in the following figures

\def\PicScale{0.28}
\hspace{-2mm}
\begin{minipage}[t]{12.2 cm}
	\begin{minipage}[t]{8em}
		\includegraphics[clip,angle=-90,scale=\PicScale]{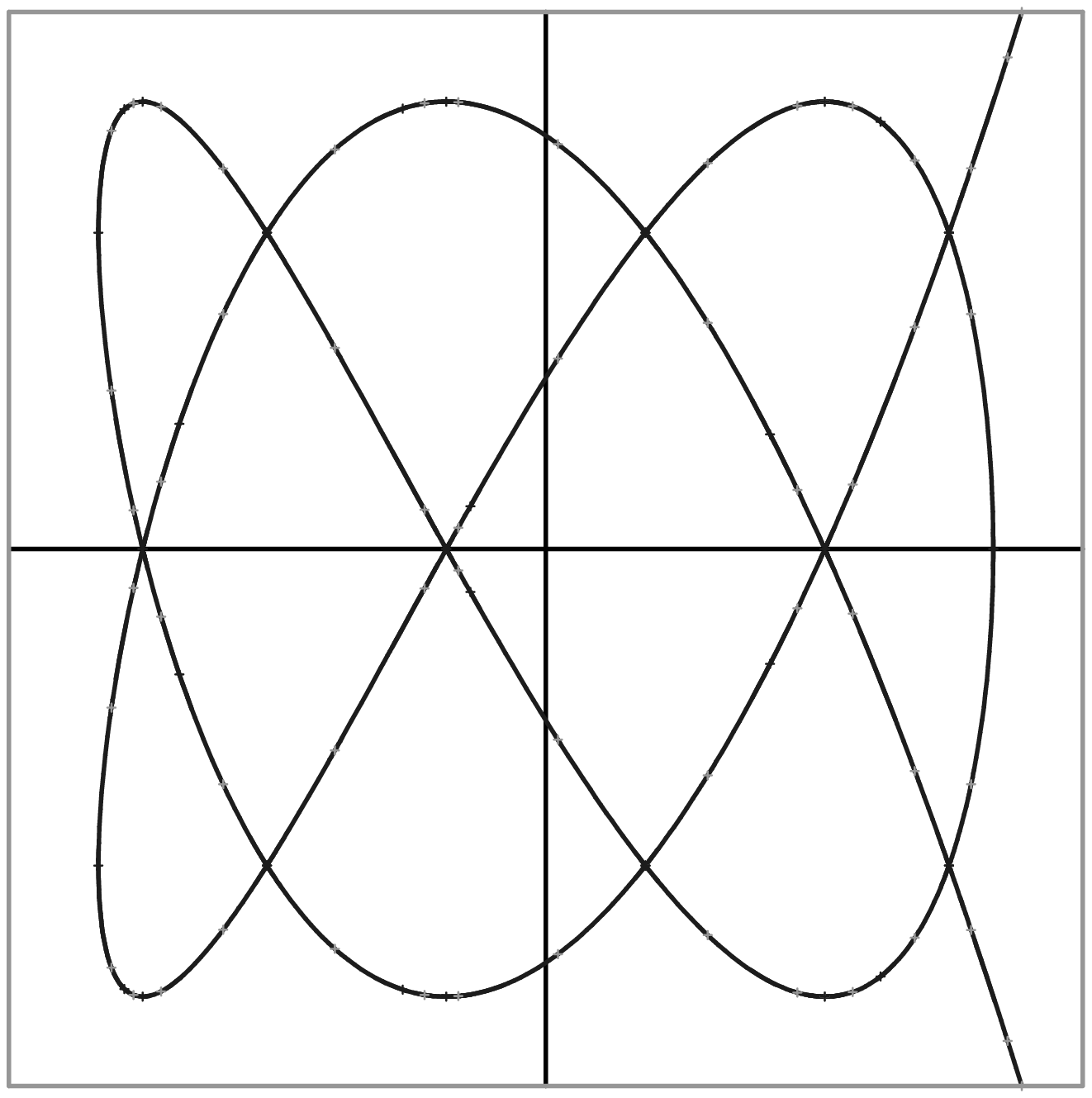}
	\end{minipage} 
	\hfill	
	\begin{minipage}[t]{8em}
		\includegraphics[clip,angle=-90,scale=\PicScale]{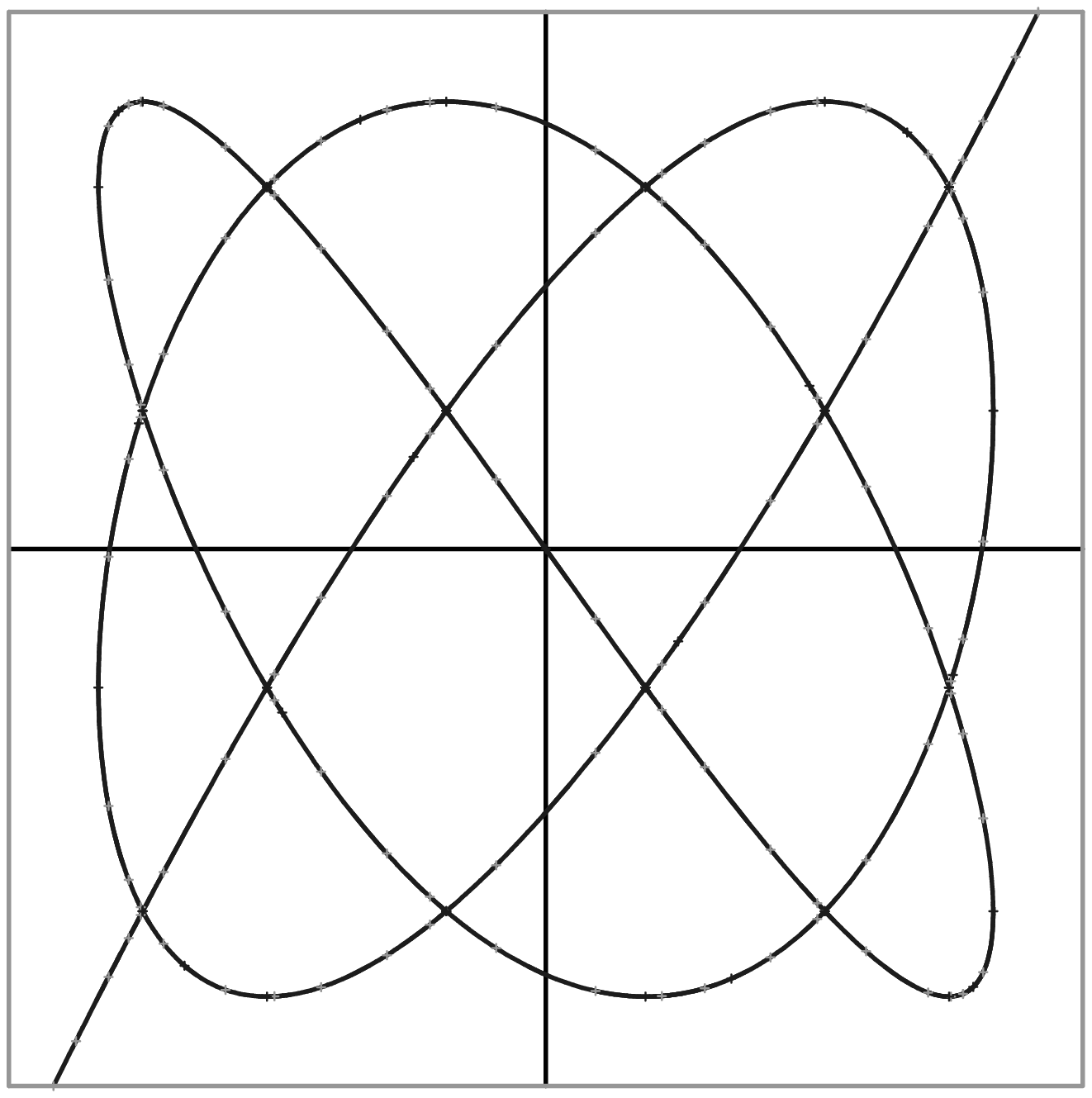}
	\end{minipage}
	\hfill
	\begin{minipage}[t]{8em}
		\includegraphics[clip,angle=-90,scale=\PicScale]{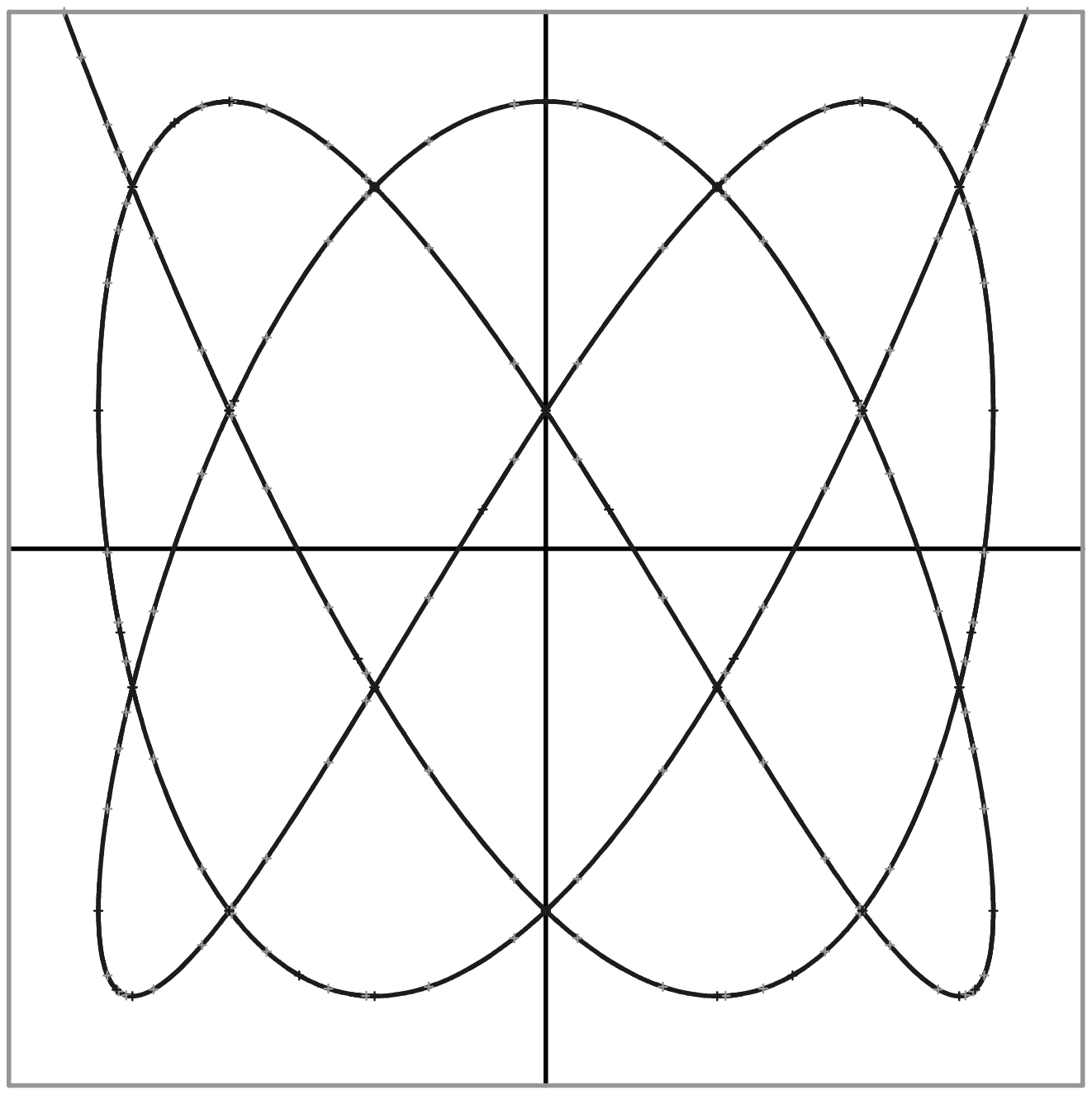}
	\end{minipage}
\end{minipage}
\def\PicScale{0.38}\\
\vspace{2mm}

The Weierstra\ss\ semigroup of $C_L$ is $\mathbb N$. Since the
$\cos \nfrac{k\pi}{q}\ (k=1,\dots,q-1)$ are algebraic numbers the zero-set of $L_{p,q}(X,Y)$ in $\mathbb A^2(K)$ defines over any 
algebraically closed base field $K$ of characteristic zero a nodal curve of type $p,q$ with $d$ distinct nodes.

Classical Lissajous curves play a role in physics in the theory of the pendulum. These curves are parametrized by the  trigonometric sine 
function. The variation of the classical curves by using Chebyshev polynomials is introduced in [Pe].
}

\headA{Elimination of nodes}\label{EliminationOfNodes}

Let $C$ be a curve of type $p,q$. Under certain assumptions the following procedure allows to eliminate singularities of $C$. 
In particular one can produce from a nodal curve $C$ further nodal curves with fewer nodes.
If the nodes of $C$ are explicitly known the new curves also have explicitly  known nodes. In this case  
one can compute their Weierstra\ss\ semigroups (see  Proposition \ref{PropLinEq}).

For relatively prime polynomials $G,H\in K[X,Y]\setminus \{0\}$ set $G_d := G+d\cdot H$ with $d\in K\setminus \{0\}$. 
Let $C:=V(G)$, $D:=V(H)$ and $C_d:=V(G_d)$ be the corresponding schemes in $\mathbb A^2(K)$ and $\Sing C$, $\Sing D$, 
$\Sing{C_d}$ their sets of singularities.
For a polynomial $\varphi\in K[X,Y]$ we denote its partial derivatives by $\varphi_X$ and $\varphi_Y$.

\prop{\label{ElimNodes}
Assume that $\Sing{C_d}$ is finite for all $d\in K\setminus \{0\}$.\\
\noindent 
a) Then for all $d\in K\setminus \{0\}$
$$
	\Sing C \cap \Sing D \subset \Sing{C_d}.
$$

\noindent
b) For almost all $d\in K\setminus \{0\}$ we have equality in a).
}

\prf{
a) $\Sing{C_d}$ is the set of common zeros of
$$
	G_d = G+dH, (G_d)_X = G_X+dH_X\ \text{and}\ (G_d)_Y = G_Y+dH_Y.
$$
For all $d\ne d'$ from $K\setminus \{0\}$ we have
$$
C_d\cap C_{d'}=C_d\cap C=C_d\cap D=C \cap D \leqno(1)
$$
since $(G_d,G_{d'})=(G_d,G)=(G_d,H)=(G,H)$, and similarly
\[
\Sing{C_d}\cap \Sing{C_{d'}}	= \Sing{C_d} \cap \Sing C \\
									=\Sing{C_d} \cap \Sing D \\
									= \Sing C \cap \Sing D. \leqno(1')
\]
Moreover these sets are finite as $G$ and $H$ are relatively prime.
In particular a) holds.\\

\noindent
{\bf Assumption:} $\Sing{C_d}\ne \Sing{C}\cap\Sing{D}$ for infinitely many $d\in K\setminus \{0\}$.

Then since $\Sing{C_d}\cap\Sing{C_{d'}} = \Sing{C}\cap\Sing{D}\ \text{for}\ d\ne d'$ the set
$$
	S:=\bigcup_{d\in K\setminus \{0\}}\Sing{C_d}
$$
is infinite. On the other hand (1') and (1) imply that
$$
	S \cap \Sing{C} = S \cap \Sing{D} = \Sing{C} \cap \Sing{D}\ \text{and}\ C_d \cap S \subset (C\cap D)\cup \Sing{C_d}\leqno(2)
$$
are finite sets for all $d \in K \setminus \{0\}$.

Let $\Gamma:=G+Z\cdot H\in K[X,Y,Z]$. Then $S$ is the set of closed points of the image of 
$V(\Gamma,\Gamma_X,\Gamma_Y) \setminus V(Z) \subset \mathbb A^3(K)$ by the orthogonal projection into the $X,Y$-plane. 
In particular $S$ is an infinite constructible set in $\mathbb A^2(K)$.

For each $d\in K$ we have
$$
	G_X\cdot H-H_X\cdot G=H\cdot(G_d)_X-H_X\cdot G_d,\ G_Y\cdot H-H_Y\cdot G=H\cdot(G_d)_Y-H_Y\cdot G_d.
$$
Therefore
$$
	\text{Each}\ (\xi,\eta)\in S\ \text{is a common zero of}\ G_X\cdot H-H_X\cdot G\ \text{and}\ G_Y\cdot H-H_Y\cdot G.\leqno(3)
$$
These polynomials cannot be relatively prime since $S$ is an infinite set. They have a common prime divisor $P$ so 
that $S\cap V(P)$ is also infinite. This is an infinite constructible subset of the curve $V(P)\subset \mathbb A^2(K)$. 
It follows that $S\cap V(P)$ is open and dense in $V(P)$.

Now we distinguish three cases:

\noindent
(i) If $P$ divides $H$, then $P$ divides $H_X\cdot G$ and $P$ divides $H_Y\cdot G$. Since $G$ and $H$ are relatively prime    
$P$ divides $H_X, P$ divides $H_Y$ and $V(P)\subset V(H,H_X,H_Y) = \Sing{D}$.
This means that $S \cap \Sing D$ contains the infinite set $S \cap V(P)$ contradicting (2).

\noindent
(ii) If $P$ divides $G$, then we obtain a contradiction against the finiteness of $S \cap \Sing{C}$.

\noindent
(iii) If $P$ does not divide $G\cdot H$ the set $V(P)\cap V(G\cdot H)$ is finite and thus the set 
$U:=(S\cap V(P))\setminus V(G\cdot H)$ is open and dense in $V(P)$. The quotient $G/H$ defines a non-vanishing regular 
function $\varphi$ on $U$. By (3) we have
$$
	(G/H)_X(\xi,\eta)=(G/H)_Y(\xi,\eta)=0\ \text{for all}\ (\xi,\eta)\in U.
$$
Therefore the differential $d\varphi$ vanishes on $U$, and $\varphi$ is constant on $U$ and not zero. Hence there exists 
$d\in K \setminus\{0\}$ with $U\subset C_d$.
But then since $U\subset S$ the set $C_d\cap S$ is infinite, again contradicting (2).

The above assumption is false, and we have shown that Sing$(C_d)=\Sing{C} \cap \Sing{D}$ for almost all $d\in K \setminus \{0\}$.
}

Let $C : G=0$ be a curve of type $p,q$, and let $\mathcal F$ denote the filtration of $K[X,Y]$ with 
$\deg[F] X=p$, $\deg[F] Y=q$ and $\deg[F] a=0$.
We set $G_d := G+d\cdot H$ where $H \in K[X,Y] \setminus \{0\}$ with $\deg[F] H <pq$ and $d \in K$. 
Then $G_d=0$ defines a curve $C_d$ of type $p,q$ and $\Sing{C_d}$ is a finite set. 
Clearly $G$ and $H$ are relatively prime. Thus the preconditions of Proposition \ref{ElimNodes} are fulfilled.

\cor{\label{CorNodeElim}
If $\Sing C \cap \Sing D$ consists of nodes of $C$, then $C_d$ is for almost all $d\in K\setminus \{0\}$ a 
nodal curve with the set of nodes $\Sing C \cap \Sing D$.
}

\prf{
From
$$
	(G_d)_{XX} = G_{XX} + dH_{XX},\ (G_d)_{XY} = G_{XY} + dH_{XY},\ (G_d)_{YY} = G_{YY} + dH_{YY}
$$
it follows for the Hesse determinants $\Hess G$ and $\Hess{{G_d}}$ of $G$ and $G_d$ that
$$
	\Hess{{G_d}}=
	\Hess G + d\left(G_{XX}H_{YY} + G_{YY}H_{XX} - 2G_{XY}H_{XY}\right) + d^2\left(H_{XX}H_{YY} - H_{XY}{}^2\right).\leqno(4)
$$
If $\Sing C \cap \Sing D$ consists of the nodes $(x_k,y_k)\ (k=1,\dots,l)$ of $C$, then
$$
	\Hess G(x_k,y_k)\ne 0\ (k=1,\dots,l).
$$ 
It follows from Proposition \ref{ElimNodes} and from (4) that for almost all $d\in K\setminus \{0\}$ we have 
$\Sing {C_d} = \Sing C \cap \Sing D$ and $\Hess{{G_d}}(x_k,y_k)\ne 0$ for $k=1,\dots,l$. Therefore the $(x_k,y_k)$ are also nodes of $C_d$.
}

Let $C$ be a nodal curve of type $p,q$ having certain nodes $(x_1,y_1),\dots,(x_l,y_l)$. We want to obtain a nodal curve 
$C_d$ such that the nodes $(x_1,y_1),\dots,(x_l,y_l)$ are preserved while all other nodes of $C$ are eliminated. 
According to the above construction we have to choose a curve $D: H=0$ having the singularities $(x_k,y_k)\ (k=1,\dots,l)$ 
and being regular at the other nodes of $C$. If we can find such an $H$ we say that $C_d$ is obtained from $C$ by {\it elimination of nodes}.

For example we can try $H=L^2$ with $L\in K[X,Y],\ 2\,\deg[F] L<pq$ where the curve $L=0$ passes 
through the nodes $(x_k,y_k)\ (k=1,\dots,l)$, but through none of the other nodes of $C$. In particular the 
curve $L=0$ can be a union of lines:
$$L=\prod_{i=1}^r(Y-m_iX+a_i)\prod_{j=1}^s(X-b_j)\ (m_i,a_i,b_j\in K)\ \text{with}\ 2(pr+qs)<pq.$$

\expl{
Elimination of the nodes of the Lissajous curve of type $3,5$
$$L_{3,5}=4Y^3-3Y-16X^5+20X^3-5X=0.$$

The figures show the real points of the curves $C_i: L_{3,5}+G_i=0\ (i=1,\dots,4)$ where
\begin{eqnarray*}
	G_1 &=& \nfrac{1}{50}\left(X - \nfrac{1}{4}\left(1 - \sqrt{5}\,\right)\right)^2\left(Y + \nfrac{1}{2}\right)^2,\\
	G_2 &=& \nfrac{1}{50}\left(Y + \nfrac{1}{2}\right)^2,\\
	G_3 &=& \nfrac{1}{50}\left(X + \nfrac{1}{4}\left(1 - \sqrt 5\,\right)\right)^2,\\
	G_4 &=& \nfrac{1}{50}.
\end{eqnarray*}

\def\PicScale{0.28}
\begin{center}
\begin{minipage}[t]{25em}
	\begin{minipage}[c]{10em}
	$C_1$: \includegraphics[clip,angle=-90,scale=\PicScale]{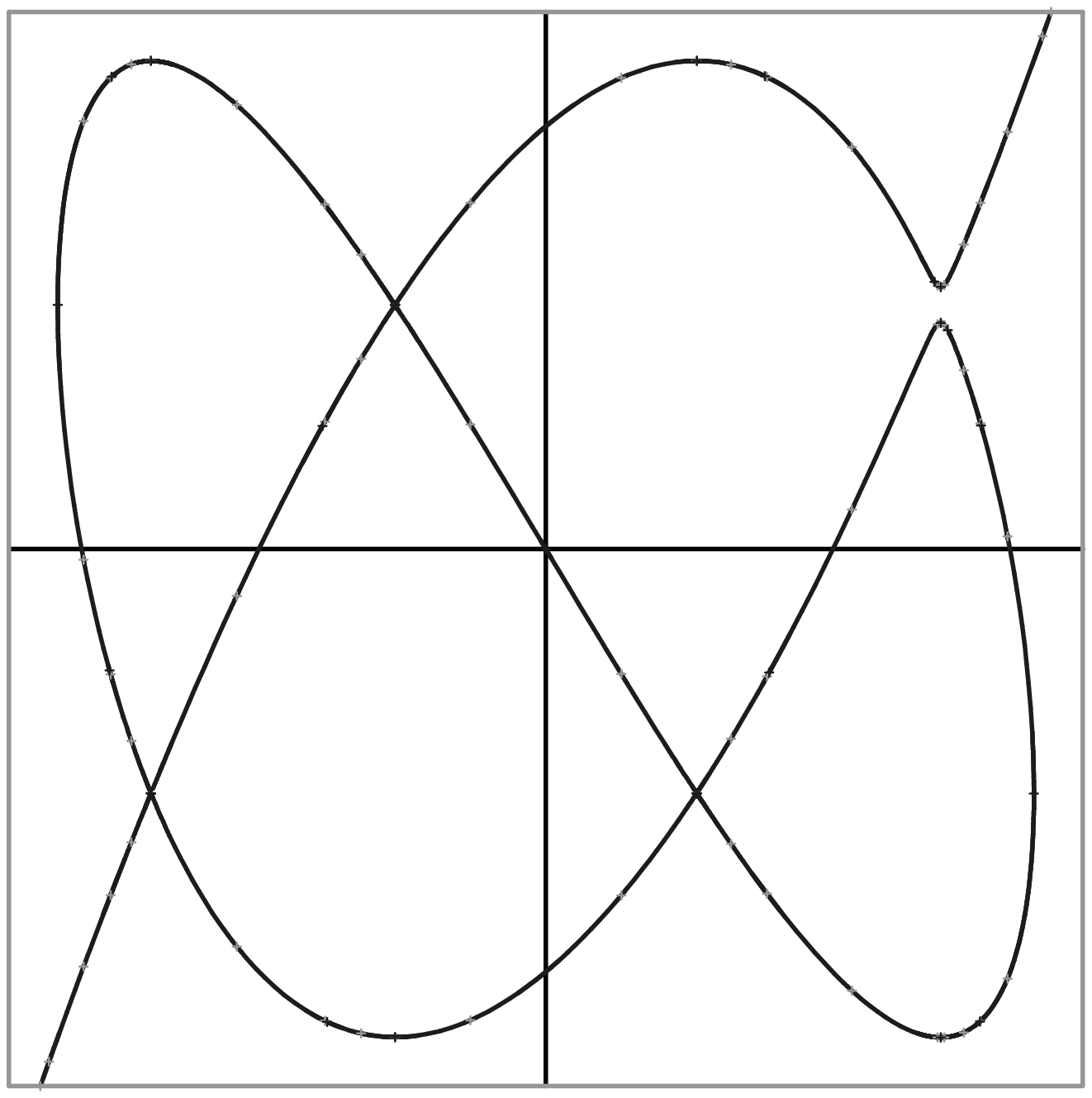}
	\end{minipage} \hfill
	\begin{minipage}[c]{10em}
	$C_2$: \includegraphics[clip,angle=-90,scale=\PicScale]{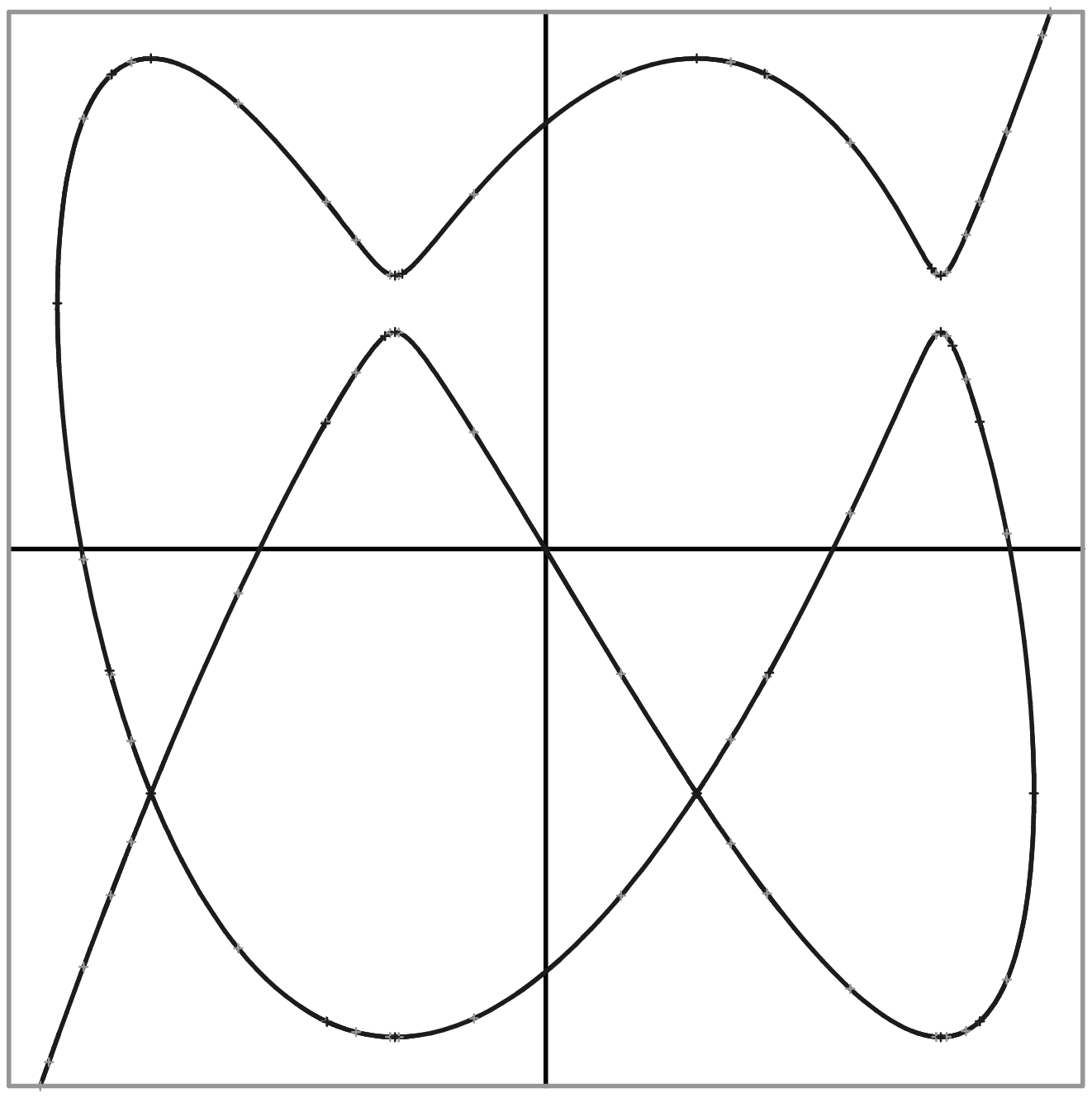}
	\end{minipage}
\end{minipage}\\

\vspace{1cm}
\begin{minipage}[t]{25em}
	\begin{minipage}[c]{10em}
	$C_3$: \includegraphics[clip,angle=-90,scale=\PicScale]{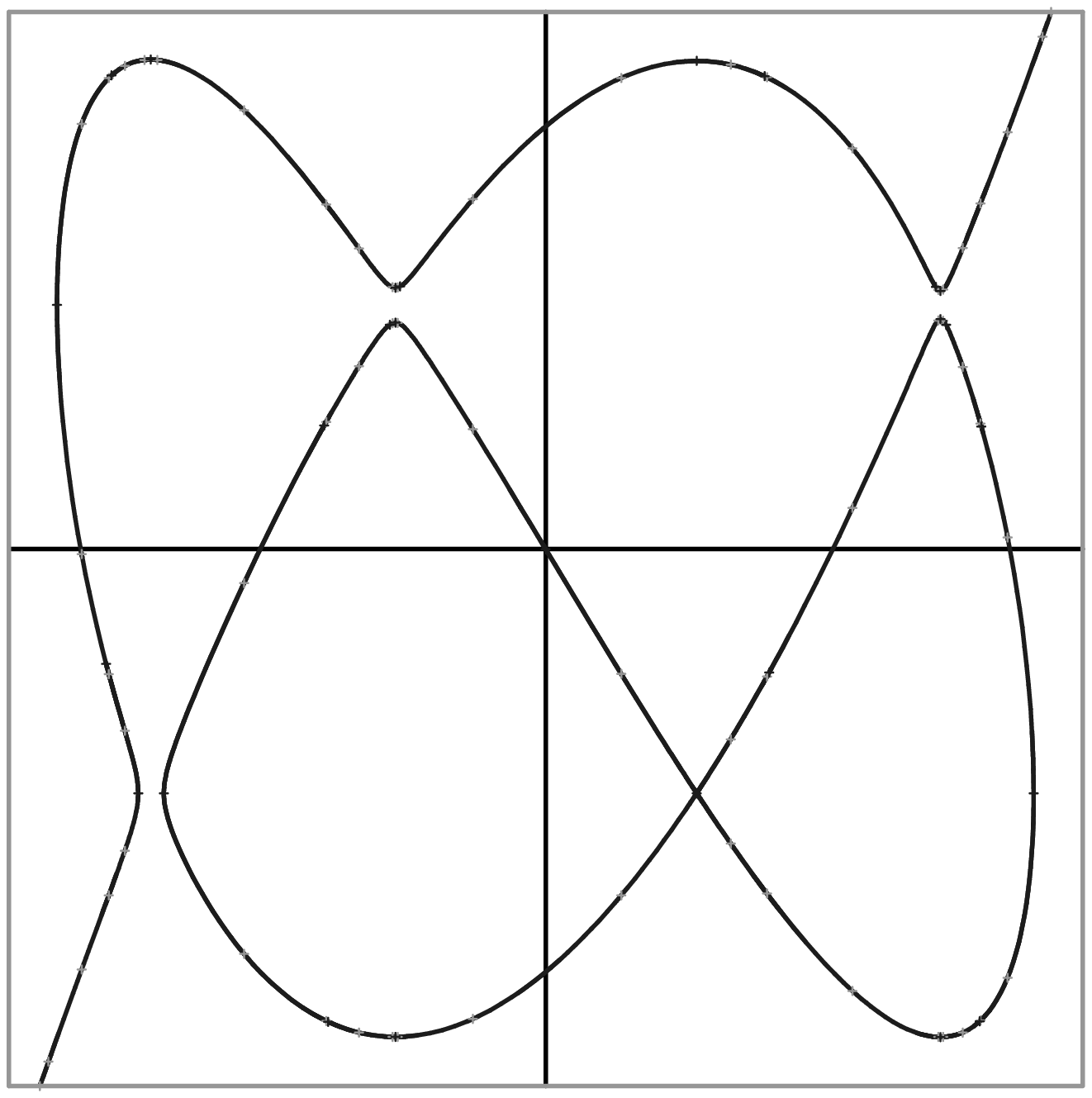}
	\end{minipage} \hfill
	\begin{minipage}[c]{10em}
	$C_4$: \includegraphics[clip,angle=-90,scale=\PicScale]{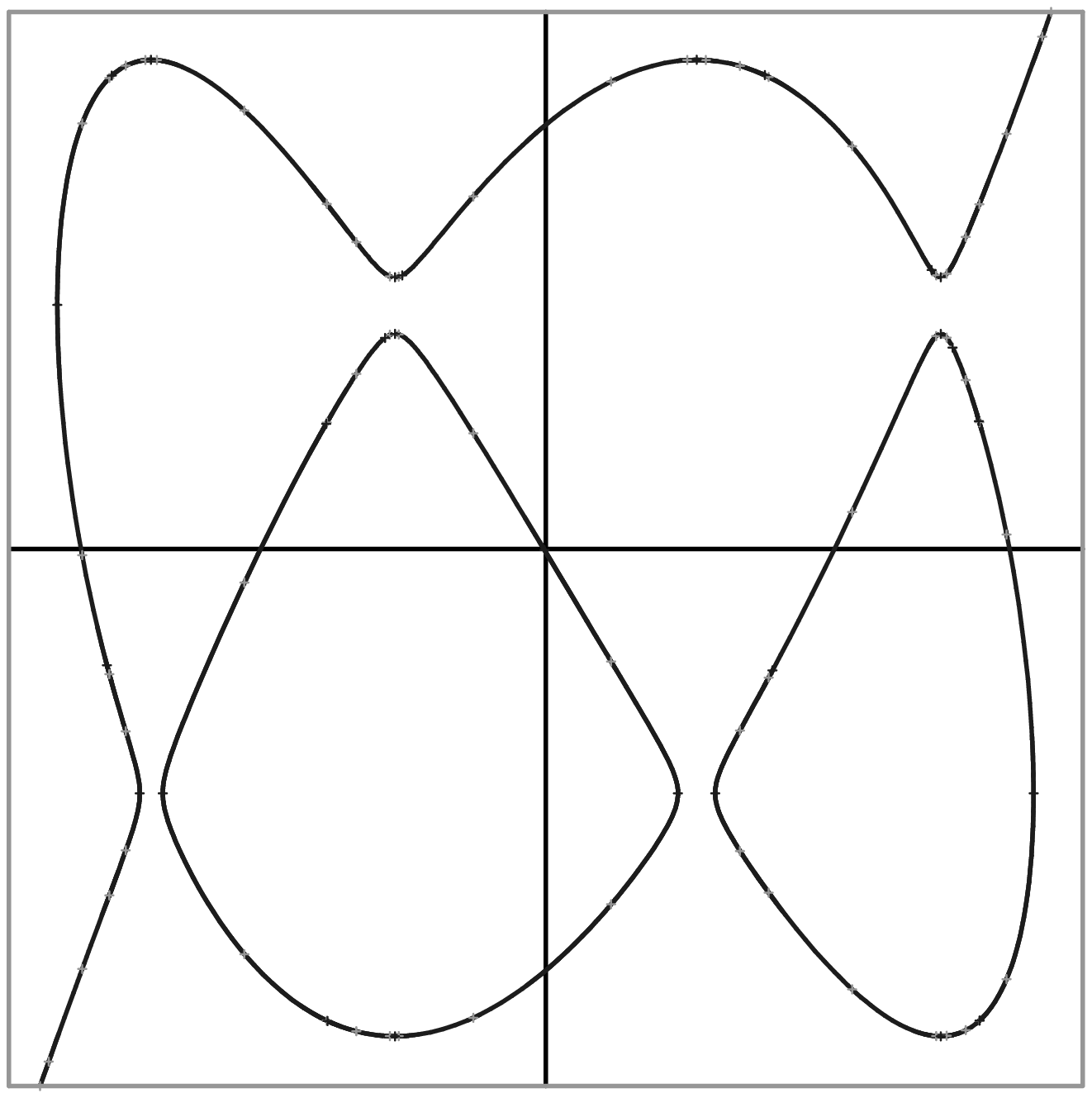}
	\end{minipage}
\end{minipage}\\
\end{center}
\def\PicScale{0.38}

\noindent
Their Weierstra\ss\ semigroups are the numerical semigroups $H$ with $3,5\in H$ except for $H=\mathbb N$.
}

If $\Gamma$ is a nodal curve of type $p,q$ whose nodes are also
nodes of another such curve $C$, then $\Gamma$ is obtained from $C$
by elimination of nodes: If $C: Y^p-X^q+H_1=0$ and $\Gamma:
Y^p-X^q+H_2=0$, take $H:=H_2-H_1$. Then $\Gamma=C_1: G+1\cdot H=0$.

\prop{\label{ElimInters}
If $\kappa_1$ and $\kappa_2$ are the sets of nodes of nodal curves of type $p,q$, then
$\kappa_1\cap \kappa_2$ is as well the set of nodes of such a curve.
} 

\prf{
Let $\kappa_i$ be the set of nodes of the nodal curve $\Gamma_i$ of type $p,q$ defined by
$F_i=Y^p-X^q+H_i\ (i=1,2)$. Set $G:=F_1, H:=H_2-H_1, D:=V(H)$. Then with the notation of \ref{ElimNodes} we have $C=\Gamma_1, C_1=\Gamma_2$ and
$\kappa_1\cap\kappa_2=\text{Sing}(C)\cap\text{Sing}(D)$ by (1'). Hence by \ref{CorNodeElim} $C_d: F_1+d\cdot H=0$ is for almost all $d\in
K\setminus\{0,1\}$ a nodal curve with the set of nodes
$\kappa_1\cap\kappa_2$.
}

Here is another class of nodal curves of type $p,q$.

\prop{\label{constrOfCurves1}
For $l_1,l_2\in \mathbb N$ with $1\le l_1 \le \nfrac{q}{2}, 1\le l_2\le \nfrac{p}{2}$ let $x_0,\dots,x_{q-l_1-1}\in K$ 
resp. $y_0,\dots,y_{p-l_2-1}\in K$ be pairwise distinct. Set
$$
	H(X):= \prod_{i=0}^{l_1-1}(X-x_i)^2\prod_{i=l_1}^{q-l_1-1}(X-x_i),\;
	G(Y):= \prod_{j=0}^{l_2-1}(Y-y_j)^2\prod_{j=l_2}^{p-l_2-1}(Y-y_j).
$$
Then $C_d:G+d\cdot H=0$ is for almost all $d\in K\setminus \{0\}$ a nodal curve of type $p,q$ with set of nodes
$$
	N:=\{(x_i,y_j)\vert i=0\dots l_1-1, j=0,\dots, l_2-1\}.
$$
}

\prf{
For all $d\in K\setminus \{0\}$ the curve $C_d$ is of type $p,q$, and $\text{Sing}(C_d)$ is finite. 
Obviously $G$ and $H$ are relatively prime. 
By Proposition \ref{ElimNodes} we have $\Sing{C_d}= \Sing C \cap \Sing D$ for almost all $d\in K\setminus \{0\}$
where $C:=V(G), D:=V(H)$. Since $(G_d)_X=d\cdot H'$ and $(G_d)_Y:=G'$ it follows that $\Sing C \cap \Sing D = N$. 
The Hesse determinant $\Hess{{G_d}}(X,Y) = d \cdot G'\!'(X) \cdot H'\!'(Y)$ does not vanish at the points of $N$, 
hence they are nodes of $C_d$.
}

The next figure shows the real points of the curve of type 5,9 with the equation
$$
	Y^2(Y-1)(Y-2)(Y-3)-X^2(X-1)^2(X-2)^2(X-3)^2(X-4)=0
$$
having 4 nodes on the $X$-axis.

\begin{center}
	\includegraphics[clip,angle=-90,scale=\PicScale]{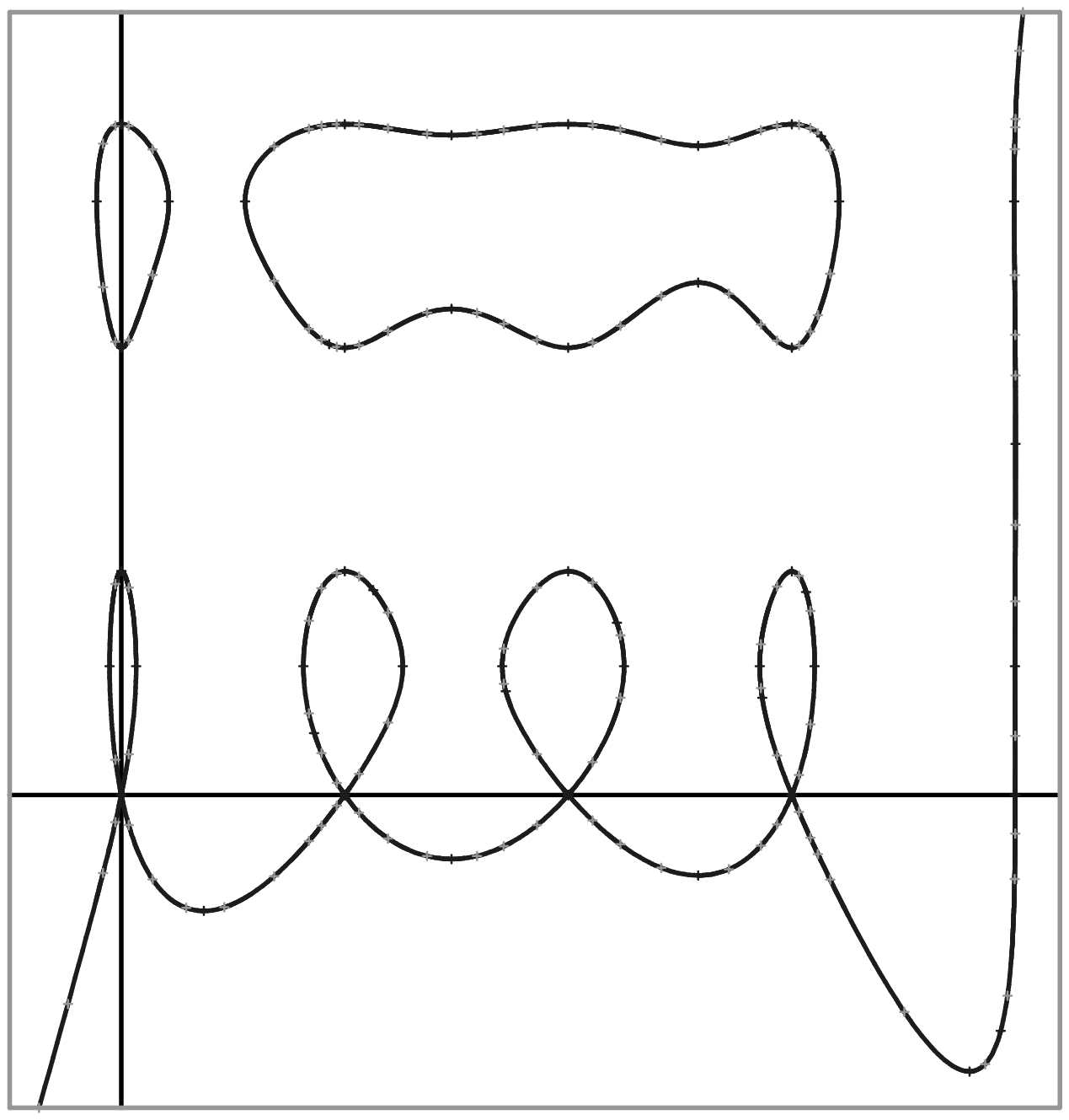}\hbox{}
\end{center}

Again for $l_1,l_2\in \mathbb N$ with $1\le l_1\le \nfrac{q}{2}, 1\le l_2 \le \nfrac{p}{2}$ let $L$ be the set of lattice points in the rectangle 
$R_{l_1-1,l_2-1}$ below the line $pX + qY = \nfrac{pq}{2}$. Set
$$
	E:=\{(\alpha,\beta)\in L\vert \alpha<l_1-1, (\alpha +1,\beta)\not\in L\}.
$$
and $L^-:=L\setminus E$. 
The following figure shows $L^-$ with $p = 43$, $q = 55$, $l_1 = 18$ and $l_2 = 15$.
\begin{center}
	\includegraphics[clip,angle=-90,scale=\PicScale]{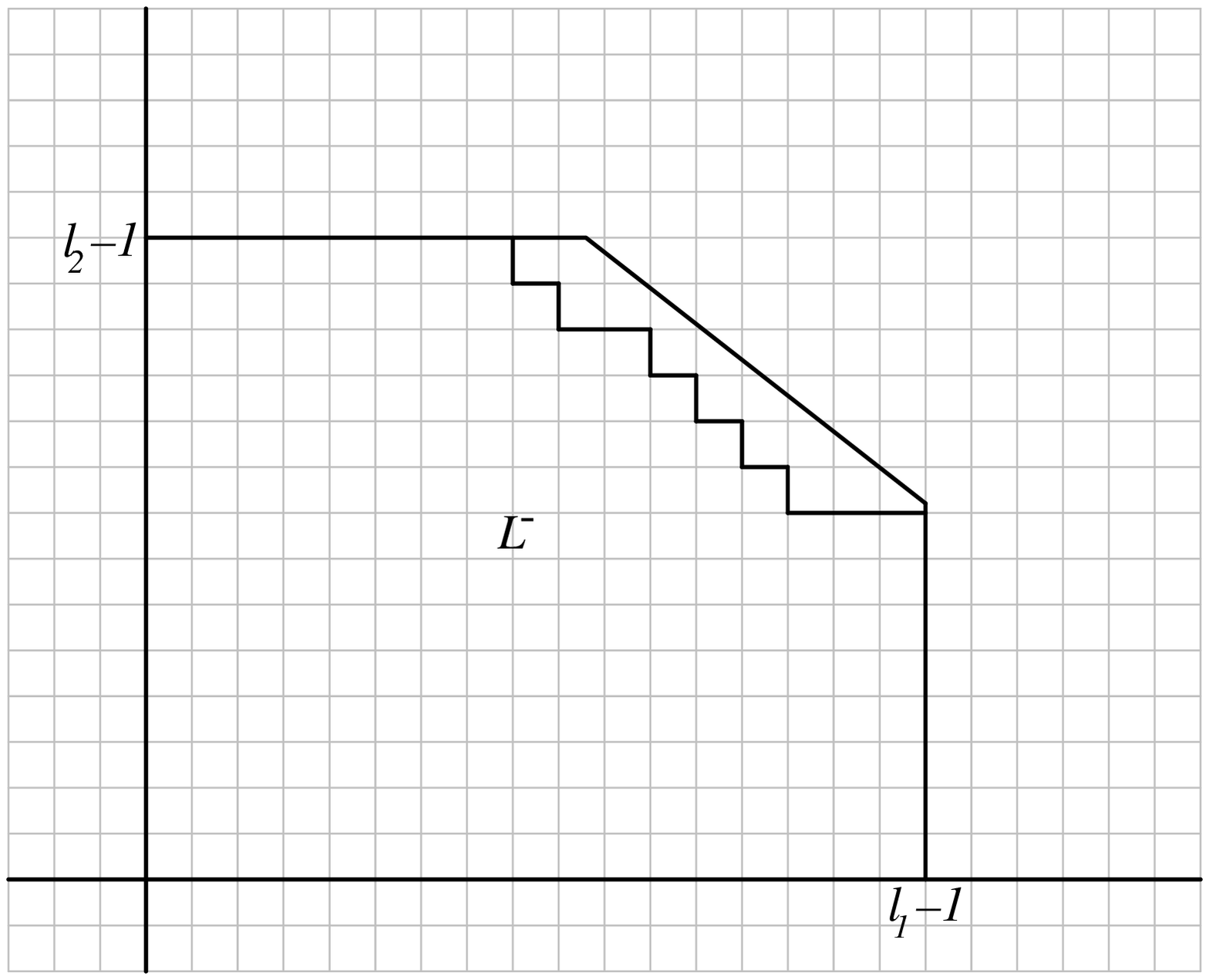}
\end{center}
\vspace{2mm}

With pairwise distinct $x_0,\dots,x_{l_1-1}\in K$ resp. 
$y_0,\dots,y_{l_2-1}\in K$ and a subset $\lambda \subset L^-$ let $\kappa(\lambda):=\{(x_m,y_n)\vert (m.n)\in \lambda\}$.

\prop{\label{constrOfCurves2}
There exists a nodal curve of type $p,q$ with the set of nodes $\kappa(\lambda)$.
}

\prf{
We apply elimination of nodes to a curve $C: F=0$ of type $p,q$ with set of nodes $N=\{(x_m,y_n)\vert m=0,\dots,l_1-1, n=0\dots,l_2-1\}$, 
see Proposition \ref{constrOfCurves1}. To this aim we construct a polynomial $H=\sum_{(i,j)\in L}u_{ij}X^iY^j\in K[X,Y]$ such that for 
$(m,n)\in R_{l_1-1,l_2-1}\cap \mathbb N^2=:R$
$$
	H(x_m,y_n)=0\ \text{if}\ (m,n)\in \lambda,\ H(x_m,y_n)\ne 0\ \text{if}\ (m,n)\not\in \lambda.
$$
We have deg$_{\mathcal F}H < \nfrac{pq}{2}$, hence $C_d: F+d\cdot H^2=0$ is for almost all $d\in K\setminus \{0\}$ 
the desired curve, see Corollary \ref{CorNodeElim}.

For $M,N \subset R$ let $\Delta_{M,N}$ denote the matrix $\left(x_m{}^iy_n{}^j\right)_{(m,n)\in M, (i,j)\in N}$ where the rows (columns) 
are ordered according to the following order of the points of R:
$$
	(0,0),(1,0),\dots,(l_1-1,0),(0,1),(1,1),\dots,(l_1-1,1),\dots,(0,l_2-1),\dots,(l_1-1,l_2-1).
$$
Since $L$ is given by a lattice path as in the second figure of Section 1 the matrix $\Delta_{L,L}$ has the form
$$
	\Delta_{L,L}= \left(\begin{array}{cccccc}
			V_{0,0}&y_0V_{0,1}&y_0{}^2V_{0,2}&\dots&y_0{}^lV_{0,l}\\
			V_{1,0}&y_1V_{1,1}&y_1{}^2V_{1,2}&\dots&y_1{}^lV_{1,l}\\
			\vdots&\vdots&\vdots&\dots&\vdots\\
			V_{l,0}&y_lV_{l,1}&y_l{}^2V_{l,2}&\dots&y_l{}^lV_{l,l}
			\end{array}\right)
$$
where $l:=l_2-1$ and the matrices $V_{i,j}$ depend only on the $x_m$. Here $V_{i,j+1}$ is obtained from $V_{i,j}$ by deleting some 
of its last columns, and $V_{i+1,j}$ by deleting the corresponding rows of $V_{i,j}$. The $V_{i,i}$ are van der Monde matrices. 
With elementary column transformations one shows that $\text{det}(\Delta_{L,L})\ne 0$.

For $(\alpha,\beta)\in E$ and $k=1,\dots,l_1-1-\alpha$ set $$L_{(\alpha,\,\beta,k)}:=(L\setminus\{(\alpha,\beta)\})\cup\{(\alpha+k,\beta)\}.$$
The matrix $\Delta_{L_{(\alpha,\,\beta,k)},L}$ is obtained from $\Delta_{L,L}$ if $V_{\beta,i}\ (i=0,\dots,l)$ is replaced by the 
matrix $\tilde V_{\beta,i}$, where in the last row of $V_{\beta,i}$ the value $x_{\alpha}$ is everywhere replaced by $x_{\alpha+k}$.
The $\tilde V_{i,i}\ (i=0,\dots,l)$ are again van der Monde matrices. Therefore $\text{det}(\Delta_{L_{(\alpha,\,\beta,k)},L})\ne 0$, too.

For any $(\alpha',\beta')\in R\setminus\lambda$ the matrix $\Delta_{\lambda\cup\{ (\alpha',\,\beta')\},L}$ is a submatrix of 
$\Delta_{L,L}$, in case $(\alpha',\beta')\in L$, or of $\Delta_{L_{(\alpha,\,\beta',k)},L}$, in case
$(\alpha',\beta')\not\in L, \alpha'=\alpha +k$ with $(\alpha,\beta')\in E$. Thus in any case
$$
	\text{rank}(\Delta_{\lambda\cup\{(\alpha',\,\beta')\},L})=\text{rank}(\Delta_{\lambda,L})+1.
$$
$\Delta_{\lambda,L}$ is the coefficient matrix of the linear system of equations for the unknowns $u_{ij}$
$$
	\sum_{(i,j)\in L}u_{ij}x_m{}^iy_n{}^j=0\qquad ((m,n)\in \lambda)\leqno(S)
$$
and $\Delta_{\lambda\cup\{(\alpha',\,\beta')\},L}$ of the system
$$
	\sum_{(i,j)\in L}u_{ij}x_m{}^iy_n{}^j=0\qquad ((m,n)\in \lambda\cup\{(\alpha',\beta')\}).\leqno(S_{(\alpha',\,\beta')})
$$
We find $u_{ij}\in K$ which solve (S) and do not solve $(S_{(\alpha',\beta')})$ for any $(\alpha',\beta')\in R\setminus\lambda$. 
Then $H=\sum_{(i,j)\in L}u_{ij}X^iY^j$ is the polynomial we were looking for.
}

\rem{\label{symSit}
If we set $E':=\{(\alpha,\beta)\in L\vert\beta<l_2-1, (\alpha,\beta+1)\not\in L\}$ and $L^=:=L\setminus E'$, then Proposition \ref{constrOfCurves2}
remains true with $L^=$ instead of $L^-$, due to the symmetry of the situation. In general $L^=$ is different from $L^-$, see the figure above.
}

\headA{Weierstra\ss\ semigroups and adjoints of nodal curves of type $p,q$}\label{WeierstrassSemigroups}

In the following let $C: F(X,Y)=0$ be a curve of type $p,q$ with place $P$ at infinity, coordinate ring $K[C]=K[x,y]$, 
projective closure $\bar C$, and let $\mathcal R$ be the normalization of $\bar C$. Let $L:=K(x,y)$ be the field of rational 
(meromorphic) functions of $C$ over $K$, let $\Omega^1_{L/K}$ be its module of differentials and $\Omega$ the K-vector space of 
holomorphic differentials on $\mathcal R$. We shall use the fact that the gaps of the Weierstra\ss\ semigroup of $C$ are the 
numbers $\ordP \omega + 1$ for the $\omega\in \Omega$.

Since $F_Y(x,y)\ne 0$ every $\omega\in \Omega^1_{L/K}$ can be written
$$
	\omega=\frac{\Phi(x,y)}{F_Y(x,y)}dx\ \text{with}\ \Phi(X,Y)\in K[X,Y].
$$

\noindent
By [G], proof of Theorem 12 we have $\omega\in \Omega$ if and only if the following two conditions are satisfied:

\noindent
(1) $\ordP \omega \ge 0$,

\noindent
(2) $\Phi(x,y)$ is contained in the conductor from $K[\widetilde C]$ to $K[C]$ where $\widetilde C$ is the normalization of $C$.

We call the curves $\Phi=0$ where $\Phi$ satisfies (2) the {\it adjoints} of $C$. If $C$ has only the nodes
$(x_1,y_1),\dots,(x_l,y_l)$ as singularities, then condition (2) is equivalent to

\noindent
(2') $\Phi(x_i,y_i)=0\ \text{for}\ i=1,\dots,l$.

For the moment we do not use the node condition. Denote by $\Omega_{\infty}$ the $K$-vector space of all $\omega\in \Omega^1_{L/K}$ 
of the form $\omega=\nfrac{\Phi(x,y)}{F_Y(x,y)}dx$ with a polynomial $\Phi$, which satisfy (1). 
Let $\gamma_1<\dots <\gamma_d$ be the gaps of $H_{pq}$: $\gamma_i=(p-1)(q-1)-1-(a_ip+b_iq)$
with a unique $(a_i,b_i)\in \mathbb N^2\ (i=1,\dots,d)$.

\lma{
The differentials
$$\omega_i:=\frac{x^{a_i}y^{b_i}}{F_Y(x,y)}dx\ (i=1,\dots,d)$$
form a basis of $\Omega_{\infty}$.
}

\prf{
We have
\begin{eqnarray*}
	\ordP{\omega_i} &=& -(a_ip+b_iq) - \ordP{F_Y(x,y)} + \ordP{dx} \\
					&=& -(a_ip+b_iq)+(p-1)q-(p+1)\\
					&=& (p-1)(q-1)-1-(a_ip+b_iq)-1=\gamma_i-1\ge 0.
\end{eqnarray*}
Thus $\omega_i\in \Omega_{\infty}\ (i=1,\dots,d)$.

Conversely let $\omega= \nfrac{\Phi(x,y)}{F_Y(x,y)}dx\in \Omega_{\infty}$.
Since $F$ is a monic polynomial in $Y$ of degree $p$ we can assume after reduction of $\Phi$ modulo $F$ that $\deg \Phi \le p-1$. Write
$$
	\Phi(x,y)=\varphi_1(x)y^{p-1}+\cdots +\varphi_p(x)\quad (\varphi_j\in K[X]).
$$
If $m\le p-1$, then different monomials of the form $x^ny^m$ have different pole orders $np+mq$ thanks to the fact that $p$ and 
$q$ are relatively prime. When $\lambda x^ay^b$ is the term of highest pole order of $\Phi$, then
$$
	\ordP \omega =-(ap+bq)+(p-1)q-(p+1)=(p-1)(q-1)-2-(ap+bq).
$$
From $\ordP \omega \ge 0$ follows that $ap+bq\le (p-1)(q-1)-2$. Then $(p-1)(q-1)-1-(ap+bq)$ is a gap of $H_{pq}$ and $a=a_i$, 
$b=b_i$ for some $i\in \{1,\dots,d\}$. For a suitable choice of $\lambda\in K$ the differential $\omega-\lambda\omega_i$ has 
greater order than $\omega$.
Recursively we obtain that $\Omega_{\infty}$ is generated by the $\omega_i$. It is clear that they are linearly independent.
}

Assume  that $C$ has the nodes $(x_1,y_1),\dots,(x_l,y_l)$ and no other singularities, and let $H$ be its Weierstra\ss\ semigroup. 
Consider for $i\in \{1,\dots,d\}$ the polynomials
$$
	\Phi_i:=X^{a_i}Y^{b_i}+u_{i+1}X^{a_{i+1}}Y^{b_{i+1}}+\dots +u_dX^{a_d}Y^{b_d}
$$
with indeterminates $u_{i+1},\dots,u_d.$

\prop{\label{PropLinEq}
$\gamma_i=(p-1)(q-1)-1-(a_ip+b_iq)$ is a gap of $H$ if and only if the system of linear equations
$$
	\Phi_i(x_j,y_j)=0\ (j=1,\dots,l)\leqno(G_i)
$$
has a solution.
}

\prf{
If $(G_i)$ has a solution $(u_{i+1},\dots,u_d)\in K^{d-i}$, then $\Phi_i$ defines an adjoint of $C$, and since $a_ip+b_iq$ 
is the pole order of $\Phi_i(x,y)$ the number $\gamma_i$ is a gap of $H$.

Conversely, if $\gamma_i$ is a gap, there must be a $\Phi_i$ of the above shape which defines an adjoint of $C$, hence $(G_i)$ 
must have a solution.
}

\cor{\label{CorLinEq}
Let $C$ and $C'$ be nodal curves of type $p,q$ with the Weierstra\ss\ semigroups $H$ resp. $H'$ where 
$H$ has genus $g$ and $H'$ has genus $g'$.  Assume that the nodes of $C'$ form a subset of the set of nodes of $C$. Then $H'\subset H$. 
If $C$ has $l$ nodes and $C'$ has $l'$ nodes, then $g'=g+(l-l')$. In particular nodal curves of type $p,q$ with the same nodes have 
the same Weierstra\ss\ semigroup.
}

\rem{
If the nodes $(x_1,y_1),\dots,(x_l,y_l)$ of a nodal curve of type $p,q$ are explicitly known, then Proposition \ref{PropLinEq} 
can be used to compute the gaps of its Weierstra\ss\ semigroup.
}

As an application of Proposition \ref{PropLinEq} we show

\prop{\label{Prop3Gen}
Let $C$ be a nodal curve of type $p,q$ with $l=l_1l_2$ nodes 
$P_{ij}:=(x_i,y_j)\ (i=0,\dots,l_1 - 1, j=0,\dots,l_2 - 1)$ where $x_0,\dots,x_{l_1-1}$ resp. $y_0,\dots,y_{l_2-1}$ 
are pairwise distinct and
$$
	l_1\le \frac{q}{2},\ l_2\le \frac{p}{2}.\leqno(\text{*})
$$

\noindent
Then 

\noindent
a) $C$ has the Weierstra\ss\ semigroup $H = H_{(l_1-1,l_2-1)}$
corresponding in the sense of Section \ref{AnIllustration} to the lattice points in the rectangle $R_{l_1-1,l_2-1}$, 
i.e.\ 
$$
	H= \left\langle p,q,pq-(l_1p+l_2q)\right\rangle.
$$ 
\noindent
b) If (*) is hurt, then no $C$ with such nodes exists.
}

\prf{
Condition (*) implies that the lattice points in the rectangle $R := R_{l_1-1,l_2-1}$ determine the numerical 
semigroup $H_{(l_1-1,l_2-1)}$ (see Example \ref{ExSemigroups}a)). If (*) is hurt, then $R$ does not define a semigroup, 
since the closedness condition of Proposition \ref{ClosedCond} is hurt.

Consider the $l\times l$-matrix $(x_i{}^my_j{}^n)$ given by the $P_{ij}$ and the lattice points $(m,n)\in R$ where we do not 
require at the moment that (*) is fulfilled. 
We order its rows and columns as in the proof of Proposition \ref{constrOfCurves2}. Then we get the matrix
$$
	\Delta:= \left ( \begin{array}{ccccc} 
			V&y_0V&y_0{}^2V&\dots&y_0{}^{l_2-1}V\\
			V&y_1V&y_1{}^2V&\dots&y_1{}^{l_2-1}V\\
			\vdots&\vdots&\vdots&\dots&\vdots\\
			V&y_{l_2-1}V&y_{l_2-1}{}^2V&\dots&y_{l_2-1}{}^{l_2-1}V
			\end{array}\right)
$$
with the van der Monde matrix
$$
	V:=\left ( \begin{array}{ccccc}
		1&x_0&x_0{}^2&\dots&x_0{}^{l_1-1}\\
		1&x_1&x_1{}^2&\dots&x_1{}^{l_1-1}\\
		\vdots&\vdots&\vdots&\dots&\vdots\\
		1&x_{l_1-1}&x_{l_1-1}{}^2&\dots&x_{l_1-1}{}^{l_1-1}
		\end{array}\right).
$$
$\Delta$ is the Kronecker product of the van der Monde matrix with respect to $y_0,\dots,y_{l_1-1}$ by $V$. 

By elementary column transformations it follows that $\text{det}(\Delta)\ne 0$ (or use the well-known formula 
$
	\det \Delta = \left(\prod_{l_1-1\ge i>j\ge 0}(x_i-x_j)\right)^{l_2}\cdot\left(\prod_{l_2-1\ge i>j\ge 0}(y_i-y_j)\right)^{l_1}
$
for Kronecker products).
To the lattice point $(l_1-1,l_2-1)\in R$ corresponds the system of linear equations
$$
	x_j{}^{l_1-1}y_k{}^{l_2-1}+\sum_{(\alpha,\,\beta)\in T}u_{\alpha,\,\beta}x_j{}^{\alpha}y_k{}^{\beta}=0\quad  (j=0,\dots,l_1-1,k=0,\dots,l_2-1)\leqno(G)
$$
where $T$ is the set of lattice points $(\alpha,\beta)$ with $(p-1)(q-1)-1-(\alpha p+\beta q) > \gamma:=(p-1)(q-1)-1-((l_1-1)p+(l_2-1)q)$. 
In particular $R \is \Nat^2 \setminus\{(l_1-1,l_2-1)\}\subset T$.

The coefficient matrix of (G) contains all columns of the Kronecker matrix $\Delta$ except for the last one, and it contains 
for the $(\alpha,\beta)\in T\setminus R$ the columns 
$s_{\alpha\beta}:=(x_j{}^{\alpha}y_k{}^{\beta})$ with $j=0,\dots,l_1-1, k=0,\dots,l_2-1$ where $\alpha<l_1-1$ or $\beta<l_2-1$. 
Let $\Delta_{\alpha\beta}$ be the matrix which is obtained from $\Delta$ by replacing its last column by $s_{\alpha\beta}$. 
We want to show that $s_{\alpha\beta}$ is linearly dependent on the other columns of $\Delta_{\alpha\beta}$. 
By symmetry we may assume that $\alpha<l_1-1$. Using $\alpha<l_1-1$ we can transform $\Delta_{\alpha\beta}$ by elementary column 
transformations into a matrix of the form
$$
	\left( \begin{array}{cccccc}
		V&0&0&\dots&0&0\\
		*&V&0&\dots&0&0\\
		\vdots&\vdots&\vdots&&\vdots&\vdots\\
		*&*&*&\dots&V&0\\
		*&*&*&\dots &*&V_{\alpha}
	\end{array}\right)
$$
where $V_{\alpha}$ is obtained from $V$ by replacing its last column by the column $(x_j{}^{\alpha})_{j=0,\dots,l_1-1}$. 
Since $\alpha<l_1-1$ we have $\det{V_{\alpha}}=0$ and hence $\det{\Delta_{\alpha\beta}}=0$.

It follows that the coefficient matrix of (G) has rank $l-1$ while $\Delta$ has rank $l$. 
Therefore (G) has no solution. By Proposition \ref{PropLinEq} $\gamma$ is an element of the Weierstra\ss\ semigroup $H$ of $C$, i.e.\ 
$H_{(l_1-1,l_2-1)}\subset H$. As $C$ has exactly $l$ nodes the number of gaps of $H$ is $\nfrac{1}{2}(p-1)(q-1)-l$, and $H_{(l_1-1,l_2-1)}=H$ follows.

If (*) is hurt we arrive at a contradiction since $R_{l_1-1,l_2-1}$ does not define a numerical semigroup, 
see Example \ref{ExSemigroups}a). Thus \ref{Prop3Gen}b) follows. If (*) is fulfilled $H$ is the semigroup generated by $p,q$ and 
$\gamma=(p-1)(q-1)-1-((l_1-1)p+(l_2-1)q)=pq-(l_1p+l_2q)$.
}

Curves $C$ as in Proposition \ref{Prop3Gen} are constructed in Proposition \ref{constrOfCurves1}. 
Observe that their Weierstra\ss\ semigroups depend only on $l_1,l_2$, 
and not on the special choice of the $x_i\ (i=0,\dots,l_1-1)$ and $y_j\ (j=0,\dots,l_2-1)$.
We can choose in particular the set of the $l_1l_2$ lattice points in the rectangle $R_{l_1-1,l_2-1}$ 
with the corners $(0,0),(l_1-1,0),(l_1-1,l_2-1),(0,l_2-1)$ as set of nodes. 
In this case we denote a curve having these nodes by $C_{l_1-1,l_2-1}$. The Weierstra\ss\ semigroup of such a curve is then the numerical semigroup 
$H_{(l_1-1,l_2-1)}$ determined in the sense of Section \ref{AnIllustration} by the same rectangle.

It is shown in [W] that all numerical semigroups which are special almost complete intersections are Weierstra\ss\ semigroups. 
In particular all numerical semigroups with 3 generators are so, since by [He] they are complete or special almost complete intersections.

\headA{Construction of Weierstra\ss\ semigroups by elimination of nodes}\label{ConstructionOfWeierstrassSemigroups}

Let $a:=\left\lfloor\frac{q}{2}\right\rfloor,b:=\left\lfloor\frac{p}{2}\right\rfloor$. 
By Example \ref{ExSemigroups}a) the subsemigroups of $H_{(a-1,b-1)}$ are 
in one-to-one correspondence with the lattice paths inside $R_{a-1,b-1}$. If $L=(P_0,\dots,P_m)$ is given by such a path, then
$$
											\alpha_0=0 < \alpha_1 < \dots < \alpha_{m-1} \le \alpha_m \le \left\lfloor\frac{q}{2}\right\rfloor-1,\ 
	\left\lfloor\frac{p}{2}\right\rfloor -1 \ge \beta_0 \ge \beta_1 > \dots > \beta_{m-1} > \beta_m = 0.
$$
Let $H_L$ denote the numerical semigroup corresponding to $L$. Then
$$
	H_L = \bigcup_{i=0}^m H_{(\alpha_i,\,\beta_i)}\leqno(1)
$$
(see the second figure in Section \ref{AnIllustration}).

Suppose there exists a nodal curve $C_L$ of type $p,q$ whose nodes are the points of $L$. Examples for this are the curves $C_{l_1-1,l_2-1}$ 
with $l_1\le a,l_2\le b$  and curves found by elimination of nodes of the $C_{l_1-1,l_2-1}$ (see Proposition \ref{constrOfCurves2}).

\thm{\label{ThmConstByElim}
$C_L$ has the Weierstra\ss\ semigroup $H_L$.
}

\prf{
The statement is true for the rectangles $R_{\alpha_i,\,\beta_i}\ (i=0,\dots,m)$ (Proposition \ref{Prop3Gen}). 
Let $H$ denote the Weierstra\ss\ semigroup of $C_L$. The nodes of $C_{\alpha_i,\,\beta_i}$ form a subset of the set of nodes of $C_L$, 
therefore $H_{(\alpha_i,\,\beta_i)}\subset H$ by \ref{CorLinEq} and hence $H_L\subset H$ by (1).

By \ref{GenusFormula} the number of gaps of $H$ is $\nfrac{1}{2}(p-1)(q-1)-l$ where $l$ is the number of nodes of $C_L$. But this is also the number 
of gaps of $H_L$, hence $H=H_L$.
}

An $L$ as in Theorem \ref{ThmConstByElim} defines a whole series of Weierstra\ss\ semigroups:

\cor{
Let $p'<q'$ be relatively prime integers with $p<p',q<q'$. Closing the gaps of $H_{p'q'}$ 
corresponding to the lattice points in $L$ gives a numerical semigroup $H'_L$. It is the Weierstra\ss\ semigroup of a 
nodal curve of type $p',q'$ with $L$ as its set of nodes.
}

\prf{
Let $a' := \left\lfloor\frac{p'}{2}\right\rfloor, b':= \left\lfloor\frac{q'}{2}\right\rfloor$. Then 
$L\subset R_{a-1,b-1}\subset R_{a'-1,b'-1}$ and $L$ defines indeed a numerical subgroup $H'_L$. Let $H$ be the Weierstra\ss\ 
polynomial defining $C_L$ and $\Phi$ that of the curve $C'$ of type $p',q'$ corresponding to $R_{a'-1,b'-1}$. Let $\mathcal F'$ 
be the degree filtration of $K[X,Y]$ with $\deg[F] X=p',\deg[F'] Y = q'$ and 
$\deg[F'] a=0$ for $a\in K$. Since $\deg[F'] H < p'q'$ and the points in $L$ are the 
singularities of $C_L$ the polynomial $H$ can be used to eliminate the nodes of $C'$ except those in $L$. We find a 
nodal curve $C'_L: \Phi+d\cdot H=0$ of type $p',q'$ which by Theorem \ref{ThmConstByElim} has the Weierstra\ss\ semigroup $H'_L$.
}

\cor{
Let $L_1$ and $L_2$ be given by lattice paths inside $R_{a-1,b-1}$. Assume that by 
eliminating the nodes of $C_{a-1,b-1}$ outside of $L_1$ resp. $L_2$ we can construct nodal curves $C_{L_1}$ and $C_{L_2}$. 
Then $H_{L_1}\cap H_{L_2}$ is a Weierstra\ss\ semigroup.
}

\prf{
This follows from \ref{ThmConstByElim} since $H_{L_1}\cap H_{L_2}=H_{L_1\cap L_2}$ and $L_1\cap L_2$ is by \ref{ElimInters} the set of nodes 
of a nodal curve of type $p,q$ obtained from $C_{a-1,b-1}$ by elimination of nodes.
}

\expl{
For $l_1,l_2\in \mathbb N$ with $1\le l_1\le a, 1\le l_2\le b$ let $L^-$ be as in Proposition \ref{constrOfCurves2}, and set 
$x_i=i\ (i=0,\dots,l_1-1),\ y_j=j\ (j=0,\dots,l_2-1)$. Then there is a nodal curve $C_{L^-}$ of type $p,q$ with the set of nodes 
$L^-$. Since $L^-$ is given by a lattice path $C_{L^-}$ has the Weierstra\ss\ semigroup $H_{L^-}$. The subsets $\lambda\subset L^-$ 
given by lattice paths are in one-to-one correspondence with the subsemigroups $H\subset H_{L^-}$ with $p,q\in H$. 
By Proposition \ref{constrOfCurves2} all these $H$ are Weierstra\ss\ semigroups, because there are curves 
$C_{\lambda}$ with $\lambda$ as their set of nodes. Similarly for $H_{L^=}$, see Remark \ref{symSit}.

For example let $l_1=a, l_2=b$ and $r<a$. The set $\lambda$ of the lattice points on or below the line $X+Y=r-1$ belongs to $L^-$, 
hence $H_{\lambda}$ is a Weierstra\ss\ semigroup.

Under the weaker assumption that $r<p-1$ it is shown in [Kn], Satz 5b) that a curve of type $p,q$ exists having only a 
singularity of multiplicity $r$ with $r$ distinct tangents such that $H_{\lambda}$ is the Weierstra\ss\ semigroup of this curve.
}

\headA{Plane models of type $p,q$ of smooth projective curves}\label{PlaneModels}

Let $\mathcal R$ be a smooth projective curve of genus $g$ and $L=K(\mathcal R)$ its field of rational 
(meromorphic) functions over $K$. For a closed point $P\in \mathcal R$ let $\widetilde {\mathcal R}:=\mathcal R\setminus P$ and
$$A:=\{f\in L \vert f\in \mathcal O_{\mathcal R,Q}\ \text{for all closed}\ Q\in \widetilde{\mathcal R}\}.$$
Then $H(P)=\{-\ordP f  \vert  f\in A\}$.

\lma{\label{NoethNorm}
a) Let $x\in A$ with $-\ordP x =: p > 0$ be given. Then $A$ is a free module over $K[x]$ of rank $p$. 
In particular $A$ is an affine $K$-algebra and $K[x]\subset A$ a Noether normalization. 
Moreover $\widetilde{\mathcal R}=\Spec A$ is a smooth affine curve and $[L:K(x)]=p$.

\noindent
b) Let $\{n_1,\dots,n_t\}$ be the minimal system of generators of $H(P)$ and let $x_i\in A$ with $-\ordP{x_i} = n_i$ $(i=1\ldots t)$ be given. 
Then $A=K[x_1,\dots,x_t]$.
}

\prf{
a) Let $h_i$ be the smallest number in $H(P)$ such that $h_i \equiv i \mod p\ (i=0,\dots,p-1)$. 
Choose $x_i \in A$ with $-\ordP{x_i} = h_i\ (i=1,\dots,p-1), x_0 := 1$.
For $f\in A$ with $-\ordP f =: \nu$ write $\nu = \mu \cdot p + h_r$ with $\mu \in \mathbb N$. Then $x^{\mu}x_r$ has pole order 
$\nu$ at $P$ and there exists $\lambda \in K$ such that $-\ordP{f-\lambda x^{\mu}x_r} < \nu$. Inductively we find elements 
$\lambda_i \in K$ such that $f - \sum\lambda_i x^{\mu} x_{r_i}$ has order 0 at $P$, hence is an element of $K$. Therefore 
$\{1,x_1,\dots,x_{p-1}\}$ is a system of generators  of $A$ over $K[x]$. As the functions have different pole 
orders mod $p$ they are linearly independent.

\noindent
b) For $f\in A$ write $-\ordP f = \sum_{i=1}^t \mu_in_i\ (\mu_i\in \mathbb N)$. Then there exists $\lambda \in K$ such that 
$f-\lambda x_1{}^{\mu_1}\cdots x_t{}^{\mu_t}$ has lower pole order at $P$ than $f$. Inductively we obtain that $f\in K[x_1,\dots,x_t]$.
}

\lma{\label{PrimitiveElem}
Let $x,y\in A$ be functions such that $p := -\ordP x$ and $q := -\ordP y$ are relatively prime. 
Then $y$ is a primitive element of $L$ over $K(x)$ whose minimal polynomial has the form
$$
	F(x,Y)=Y^p+f_1(x)Y^{p-1}+\cdots +f_p(x)
$$
where $f_1,\dots,f_p\in K[X],\ \deg{f_j} \le \left\lfloor\nfrac{jq\,}{p}\right\rfloor\ (j=1,\dots,p-1),\ \deg{f_p} = q$. 
Hence $F(X,Y)$ is a Weierstra\ss\ polynomial of type $p,q$.
}

\prf{
As $K[x]$ is integrally closed in $K(x)$ the minimal polynomial
$$
	F(x,Y)=Y^s+f_1(x)Y^{s-1}+\cdots +f_s(x)
$$
of $y$ has coefficients $f_i(x)\in K[x]$. By Lemma 6.1a) we have $s\le p$. Since $F(x,y)=0$ at least two terms of $F(x,y)$ 
have the same order at $P$. Assume
$$
	\text{deg}(f_j)\cdot p+q\cdot(s-j)=\text{deg}(f_l)\cdot p+q\cdot (s-l)
$$
with $j,l\in \{0,\dots,s\},\ l>j$, where we have set $f_0=1$. Then $q\cdot(l-j)\equiv 0\ \mod p$. As $p$ and $q$ 
are relatively prime this is only possible when $s=p,\ l=p$ and $j=0$. From $[L:K(x)]=p$ we conclude that $y$ is a 
primitive element of $L$ over $K(x)$.

The orders of the terms $f_1(x)y^{p-1},\dots,f_1(x)$ are pairwise distinct, and likewise $y^p$ and $f_j(x)y^{p-j}$ 
have different order for $j=1,\dots,p-1$. We obtain that $\ordP{y^p}=\ordP{f_p(x)}$ and that $y^p$ and $f_p(x)$ 
are the terms of lowest order. It follows that $\deg{f_p}=q$ and
$$
	-\ordP{f_j(x)y^{p-j}}=\deg{f_j}\cdot p+(p-j)\cdot q<p\cdot q
$$
hence $\deg{f_j} \le \left\lfloor\nfrac{jq\,}{p}\right\rfloor$ for $j=1,\dots,p-1$.
}

Lemma \ref{PrimitiveElem} shows that $\mathcal R$ has a {\it plane model} $C$ of type $p,q$, i.e.\ 
that $\mathcal R$ is the normalization of the projective closure of $C$, a fact which is covered by the older 
literature (see e.g.\ [Ha] or [HL]). But we want even a plane model which is a nodal curve of type $p,q$.

\prop{\label{SpaceModel}
Let $A=K[x_1,\dots,x_t]$ and let $p,q\in H(P)$ be relatively prime numbers such that
$$
	-\ordP{x_i} < p < q \ (i=1,\dots,t).
$$
Then there exist elements $u_1,u_2,u_3\in A$ with
$$
	-\ordP{u_1}=p,-\ordP{u_2}=q,-\ordP{u_3}<p
$$
such that $A=K[u_1,u_2,u_3]$. In particular $C=\Spec A$ can be embedded into $\mathbb A^3(K)$.
}

\prf{
Choose $u_1\in A$ with $-\ordP {u_1} = p$. Then by Lemma \ref{NoethNorm}a) $K[u_1]$ is a Noether normalization of $A$. 
Choose $y\in A$ with $-\ordP y = q$ and set $u_2 := y + \sum_{i=1}^ta_ix_i
\ (a_i\in K)$. Then $-\ordP {u_2} = q$ and $L = K(u_1,u_2)$ by Lemma \ref{PrimitiveElem}.

If $\mathfrak m \in \Max A$ is unramified over $K[u_1]$, then $K[u_1]$ contains a local parameter 
for $\mathfrak m$. Now let $\mathfrak m$ be ramified and $\tau\in A$ a local parameter for $\mathfrak m$. Write
\begin{eqnarray*}
	x_i &=& \alpha_{i0}+\alpha_{i1}\tau + \dots\quad (i=1,\dots,t)\\
	y	&=& \beta_0+\beta_1\tau +\dots
\end{eqnarray*}
$(\alpha_{i0},\alpha_{i1},\beta_0,\beta_1\in K)$, where the dots represent terms of higher order in $\tau$. 
At least one $\alpha_{i1}\ne 0$, since the $x_i$ generate $A$ as a $K$-algebra. We see that 
$u_2-(\beta_0+\sum_{i=1}^t\alpha_{i0}a_i)$ is a local parameter of $\mathfrak M$ if $(a_1,\dots,a_t)$ is 
chosen in the complement of a hyperplane in $\mathbb A^t(K)$. As $A$ is ramified over $K[u_1]$ only at 
finitely many maximal ideals we can arrange by choosing $(a_1,\dots,a_t)$ in the complement of finitely 
many hyperplanes that $K[u_1,u_2]$ contains for every $\mathfrak m \in \Max A$ a local parameter.
Set $B:=K[u_1,u_2]$. 

We are looking now for a $u_3\in A$ of the form
$$
	u_3=\sum_{i=1}^tb_ix_i,\ (b_1,\dots,b_t)\in K^t\setminus \{0\}
$$
such that different maximal ideals of $A$ are lying over different maximal ideals of $K[u_1,u_2,u_3]$.

Consider points $\mathfrak M_1\ne\mathfrak M_2$ on $\Max A$. Then
$$
	\mathfrak M_j=(x_1-\lambda_{j1},\dots,x_t-\lambda_{jt})\ (\lambda_{ji}\in K,\ j=1,2,\ i=1,\dots,t)
$$
with $(\delta_1,\dots,\delta_t):=(\lambda_{11}-\lambda_{21},\dots,\lambda_{1t}-\lambda_{2t})\ne 0$. 
If already $\mathfrak M_1\cap B\ne \mathfrak M_2\cap B$, then $\mathfrak M_1\cap K[u_1,u_2,u_3] \not= \mathfrak M_2\cap K[u_1,u_2,u_3]$ 
for arbitrarily chosen $(b_1,\dots,b_t)\in K^t$. Otherwise $\mathfrak M_1\cap B=\mathfrak M_2\cap B$ is a singular point of $\Max B$. 
For any choice of $(b_1,\dots,b_t)$ we have $-\ordP{u_3} < p$ and
$$
	\mathfrak M_j\cap K[u_1,u_2,u_3]=(u_1-\bar u_1,u_2-\bar u_2,u_3-\sum_{j=1}^t\lambda_{ij}b_i)\ (j=1,2)
$$
where $\bar u_1,\bar u_2\in K$ do not dependent on $j$. Then
$$
	\mathfrak M_1\cap K[u_1,u_2,u_3]\ne \mathfrak M_2\cap K[u_1,u_2,u_3]\ \text{if and only if}\ \sum_{i=1}^t\delta_ib_i\ne 0,
$$
i.e. if and only if $(b_1,\dots,b_t)$ is chosen outside of a hyperplane in $K^t$. As \Max B has only finitely many 
singularities and over each of them there are only finitely many $\mathfrak M\in \Max A$ we can choose $u_3$ as desired.

Set $D:=K[u_1,u_2,u_3]$. Then for any $\mathfrak m\in \Max D$ there is exactly one $\mathfrak M\in \Max A$ 
lying over $\mathfrak m$. Therefore $A_{\mathfrak m}=A_{\mathfrak M}$ and $A_{\mathfrak m}$ is a finite $D_{\mathfrak m}$-module. 
Moreover any $\mathfrak M\in \Max A$ contains a local parameter $\tau\in \mathfrak m$. Then $A_{\mathfrak m}=D_{\mathfrak m}$ 
and consequently $A=D$.
}

\thm{\label{NodalModels}
Let $p<q$ be relatively prime numbers in $H(P)$ which are greater than the elements of the minimal system of 
generators of $H(P)$. Then $\mathcal R$ has a plane model of type $p,q$ with the place $P$ at infinity which 
has at most nodes as singularities. In particular every Weierstra\ss\ semigroup is the Weierstra\ss\ semigroup 
of a nodal curve of type $p,q$ for properly chosen $p,q$.
}

\prf{
Let $\{n_1,\dots,n_t\}$ be the minimal system of generators of $H(P)$ and let $x_i\in A$ with \linebreak $-\ordP{x_i} = n_i$. 
By Lemma \ref{NoethNorm}b) we have $A=K[x_1,\dots,x_t]$ and with $u_1,u_2,u_3$ as in Proposition \ref{SpaceModel} 
we obtain $A=K[u_1,u_2,u_3]$.

It is well-known that the generic projection of $\Spec A$ from 3-space into the plane is birationally equivalent to $\Spec A$ 
and has at most nodes as singularities. Its coordinate ring has the form $K[u_1-a_1u_3,u_2-a_2u_3]$ where $(a_1,a_2)$ can be 
chosen in the complement of an algebraic curve in $\mathbb A^2(K)$ (see [A], section II of its appendix where also the case 
characteristic $K>0$ is settled). Set $u:=u_1-a_1u_3, v:=u_2-a_2u_3$. Then $-\ordP u=p$, $-\ordP v=q$. 
By Lemma \ref{PrimitiveElem} the minimal polynomial of $v$ over $K[u]$ is a Weierstra\ss\ polynomial of type $p,q$. 
Hence $\Spec{K[u,v]}$ is a plane model of $\mathcal R$ of type $p,q$ having at most nodes as singularities.
}


\end{document}